\documentclass[UTF-8,reqno]{amsart}
\usepackage{enumerate}
\usepackage{mhequ}
\usepackage{amssymb,url,color, booktabs,nccmath}
\usepackage{xcolor}
\usepackage[left=3.2cm,right=3.2cm,top=4cm,bottom=4cm]{geometry}
\usepackage{mathrsfs}
\usepackage{enumitem,dsfont}

\usepackage{color}
\usepackage[colorlinks=true]{hyperref}
\hypersetup{
    linkcolor=blue,          
    citecolor=red,        
    filecolor=blue,      
    urlcolor=cyan}

\definecolor{darkergreen}{rgb}{0.0, 0.5, 0.0}

\newtheorem{theorem}{Theorem}[section]
\newtheorem{lemma}[theorem]{Lemma}
\newtheorem{proposition}[theorem]{Proposition}

\newtheorem{corollary}[theorem]{Corollary}

\newtheorem{definition}{Definition}
\theoremstyle{definition}
\newtheorem{remark}{Remark}

\allowdisplaybreaks

\usepackage{esint}
\usepackage{cite}
\usepackage{verbatim}
\usepackage{times}
\usepackage{color}
\usepackage{bm}
\usepackage{fancyhdr}
\numberwithin{equation}{section}

\begin{document}
\title{Regularizing the Navier-Stokes Equations via Additive Type Noise plus Small Linear Multiplicative Noise}
\author{Weiquan Chen}
\address[Weiquan Chen]{Academy of Mathematics and Systems Science,
Chinese Academy of Sciences, Beijing 100190, China}
\email{chenweiquan@amss.ac.cn}

\begin{abstract}
    We derive a series of ``regularization by noise in high probability" results for incompressible Navier-Stokes equations $(d\geq3)$, via a new mechanism based on fixed point argument: For initial values in critical Besov space $\dot{B}_{q,\infty}^{-1+d/q}~ (1\leq q<\infty)$, we can expect unique global solutions with $C_{w^*,loc}\big([0,\infty);\dot{B}_{q,\infty}^{-1+d/q}\big)$-paths in high probability if the equations are perturbed by well-constructed additive type noise; for initial values in the ``end-point space" $\dot{B}_{\infty,\infty}^{-1}$ where the equations are locally ill-posed in deterministic setting \cite{BP08}, we can still expect unique solutions on $[0,T]$ (with ``modified" $*$-weak continuous paths) in high probability for large enough $T>0$, if the equations are perturbed by a well-constructed additive type noise plus a ``small" linear multiplicative noise. This means that the regularizing effect can be strengthen by small additional linear multiplicative noise in the sense that the Picard iteration remains effective for values in the largest critical space $\dot{B}_{\infty,\infty}^{-1}$. This work has nothing to do with AI models.
    
\end{abstract}\
\maketitle

\section{Introduction}\label{sec. introd.}
\setcounter{section}{1}

\par The incompressible Navier-Stokes equations are well-known for its wide applications in Fluid Mechanics and the longstanding regularity and uniqueness issues in mathematics. Over the years, we have series of breakthroughs coming up on the negative side, concerning non-uniqueness of weak solutions with the absent of external force term (see for example \cite{BCV18,BV19,BVannals19,CL22}) and non-uniqueness of Leray-Hopf weak solutions with specific non-zero external force (see \cite{ABC21,ABC23}). See also the more recent breakthrough \cite{CP25} on non-unique smooth solutions arising from specific (large) initial values in the critical space $BMO^{-1}$. These work suggest that the equations may ``break down" at least under certain theoretical circumstance. Over the decades, however, people have been striving for the positive direction, trying to figure out circumstances in which the equations have unique solutions. Among these, the most famous part may probably be the theory of small solutions in the following critical spaces:
\begin{align}
    \dot{H}^{d/2-1}\hookrightarrow L^d\hookrightarrow \dot{B}_{q,\infty}^{-1+d/q}~(1\leq q<\infty)\hookrightarrow BMO^{-1}\hookrightarrow \dot{B}_{\infty,\infty}^{-1}~.\nonumber
\end{align}
For global well-posedness results on small initial values in the above spaces (except $\dot{B}_{\infty,\infty}^{-1}$), we refer to \cite{FK64,Che92},\cite{Kat84},\cite{CMP93} and \cite{KT01} respectively. The $BMO^{-1}$ space is the largest critical space where existence results are available and is commonly considered as the largest one where one can safely apply the Picard iteration or  fixed point argument to construct solutions \cite{KT01}. The end point space $\dot{B}_{\infty,\infty}^{-1}$ on the right-hand-side is the largest possible critical space in the sense that all scaling invariant spaces for the Navier-Stokes equations are embedded in it \cite{Can03,BCD11}. However, it turns out to be ``too large" for the equations: in \cite{BP08}, J. Bourgain and N. Pavlovi\'c show the famous result of norm inflation phenomenon, namely initial values arbitrarily small in $\dot{B}_{\infty,\infty}^{-1}$ can produce solutions that become arbitrarily large in the $\dot{B}_{\infty,\infty}^{-1}$-norm by arbitrarily short time. So this means that the equations are locally (in time) ill-posed in this space in a very strong sense.


\par A challenging and active direction that combine this well-posedness issue with probability theory is the so-called regularization by noise effect. It was motivated by a classical result \cite{Ver81} in stochastic analysis that established well-posedness of the stochastic differential equation:
\begin{align}
    {\rm d}X(t)=b(t,X(t)){\rm d}t+\sigma(t,X(t)){\rm dB}_t\nonumber
\end{align}
when $b$ is just a bounded function and $\sigma$ satisfies certain conditions which include the simplest and probably most interesting case $\sigma\equiv const.\neq0$. Such noise term is independent of the solution itself and is said to be of ``additive type". In this direction, people try to ``restore" the well-posedness of Navier-Stokes equations by asking whether one can introduce a reasonable noise term so that, in appropriate sense, the resulting stochastic equations have unique solutions for general initial values. Here, ``in appropriate sense" often means ``almost surely" (i.e. with probability one) or ``with high probability" (i.e. with probability being arbitrarily close to one if the noise parameters are large enough). Regularization by additive noise for Navier-Stokes remains open and challenging, but we should mention that there are interesting results like ``almost sure absent of spacial singularity for each fixed time" \cite{FR02}, ``continuity of solution map on the initial values via Markov selection" \cite{DPD03,FR08}, etc.  We also mention here some famous results on regularization effect by other noise forms: see \cite{GHV14,RZZ14} for linear multiplicative noise; \cite{FL21,Agresti2026} for transport noise and \cite{HLL2025} for a non-linear and non-local type of noise. 

\par {\bf Our work:} In this paper, we construct a specific type of additive noise (by which we mean that the noise term is independent of the solutions themselves) so that the system is globally well-posed with high probability in the critical Besov space $\dot{B}_{q,\infty}^{-1+d/q}$ when $1\leq q<\infty$. This answers the open question of regularization by additive noise for the first time, although the general additive noise case still remains open. As a further result, we show that if a small linear multiplicative noise is added additionally into the equations, then the regularizing effect can be strengthen in the sense that: for general data in the largest critical space $\dot{B}_{\infty,\infty}^{-1}$, one can establish well-posedness theorem on any large enough time interval $[0,T]$ with high probability. This can be viewed as a restoration of the local ill-posedness result \cite{BP08} on the deterministic side.

\par {\bf Main idea:} Our strategy is rather simple compared to the mentioned literature above. The key is just that we need to find suitable solution spaces so that, when the solution is viewed as the fixed point of a bilinear map, the map should be continuous on the space and, at the same time, we should be able to control the stochastic Stokes flow in the space via the additive part of the noise that we impose into the Navier-Stokes equations. For initial values in $\dot{B}_{q,\infty}^{-1+d/q}$ when $1\leq q<\infty$, this is relatively simple as we just need to look at solution space like $L^p_t\dot{B}_{q,\infty}^{-1+2/p+d/q}$ (which is the $L^p$ version of the commonly known mild solution space $L^\infty_t\dot{B}_{q,\infty}^{-1+d/q}$). The solutions have $*$-weak continuous paths in space $\dot{B}_{q,\infty}^{-1+d/q}$ as supposed. For initial values in $\dot{B}_{\infty,\infty}^{-1}$, we need to impose a further linear multiplicative noise (with arbitrarily small but non-zero intensity) into the system, so that a corrector is impose into the bilinear map via Duhamel principle and the map become continuous on the modified solution space like $L^p_t\dot{B}_{\infty,\infty}^{-1+2/p+\delta}$ ($\delta>0$ is small). Consequently, the solution $u$ have ``modified" $*$-weak continuous paths in the sense that $u\in C_{w^*}\big([0,T];\dot{B}_{\infty,\infty}^{-1}\big)+ C_{w^*}\big([0,T];\dot{B}_{\infty,\infty}^{-1+\delta}\big)$.

\
\subsection{Main results}
\par We study Cauchy problem of the $d$-dimensional 
incompressible stochastic Navier-Stokes equations posed on whole space $\mathbb{R}^d$ ($d\geq3$):
\\
\begin{equation}\label{sto. NS}
	\left\{\begin{aligned}
		&{\rm d}u(t)=\Big(-{\rm div}\big(u(t)\otimes u(t)\big)+\nu \Delta u(t)-\nabla p(t)\Big){\rm d}t+{\rm d}\mathcal N(t)\\
		&{\rm div}\ u(t)=0\\
		&u(0)=u_0
	\end{aligned}\right. 
\end{equation}
\\
where initial value $u_0\in \dot{B}_{q,\infty}^{-1+d/q}$ ($1\leq q\leq\infty$) is divergence-free in weak sense. The system is driven by a stochastic noise term $\mathcal N$ which will be specified in different situations. In unified form, our noise term consists of an additive part and a ``small" linear multiplicative noise:
\begin{equation}
    \begin{aligned}\label{general noise form}
	\mathcal N(t)
    &=\varepsilon_{lin}\int_0^t\big(u(s)-g(s)\big){\rm dB}_s+\underbrace{\sum_{j}\int_0^tG_j(s){\rm dB}_s^j}_{=:~\mathcal N^{add}(t)}~,\\
    g(t)
    &=e^{\nu t\Delta}u_0+\int_0^te^{\nu(t-s)\Delta}{\rm d}\mathcal N^{add }(s)~,
    \end{aligned}
\end{equation}
where $\varepsilon_{lin}\in (-\nu^{1/2},\nu^{1/2})$ is constant and $\big\{{\rm B},{\rm B}^j\big\}$ is an independent family of one dimensional standard Brownian motions given on a probability basis $\left(\Omega, \mathscr{F}, \{\mathscr{F}_t\}_{t\geq0}, \mathbb{P}\right)$. The coefficients $G_j: \Omega\times[0,+\infty)\longrightarrow \dot{B}_{q,\infty}^{-1+d/q}$ would be constructed to be divergence-free and adaptive with respect to $\{\mathscr{F}_t\}$. Stochastic integrals would always be understood in the sense of It\^o integral in this paper. We shall first define the solutions via Duhamel principal. So let's denote ($\varepsilon\in\mathbb R$)
\begin{align}\label{def. D_s,t}
    \mathbf D_{s,t}(\varepsilon,{\rm B}):=\exp{\left[\varepsilon({\rm B}_t-{\rm B}_s)-\frac{\varepsilon^2}{2}(t-s)\right]}~,\quad 0\leq s\leq t<\infty~.
\end{align}
By subtracting the equation of $g$ from (\ref{sto. NS}), which is the stochastic Stokes equations
\begin{equation}\label{sto. Stokes}
	\left\{\begin{aligned}
		&{\rm d}g(t)=\nu\Delta g(t){\rm d}t+{\rm d}\mathcal N^{add}(t)\\
		&{\rm div}\ g(t)=0\\
		&g(0)=u_0~,
	\end{aligned}\right. 
\end{equation}
we see that $v:=u-g$ solves
\begin{equation}
	\left\{\begin{aligned}
		&{\rm d}v(t)=\Big(-{\rm div}\big(u(t)\otimes u(t)\big)+\nu \Delta v(t)-\nabla p(t)\Big){\rm d}t+\varepsilon_{lin} v(t){\rm dB}_t\\
		&{\rm div}\ v(t)=0\\
		&v(0)=0~,
	\end{aligned}\right.\nonumber
\end{equation}
i.e.
\begin{align}
    v(t)=-\int_0^t\mathbf D_{s,t}(\varepsilon_{lin},{\rm B})e^{\nu(t-s)\Delta}\mathcal{P}{\rm div}\big(u(s)\otimes u(s)\big){\rm d}s~.\nonumber
\end{align}
The Leray projector $\mathcal{P}={\bf Id}+\nabla(-\Delta)^{-1}{\rm div}$ is defined via Fourier multiplier. This motivates the following definition.

\begin{definition}\label{def. solu. path}
Let $p,q\in[1,\infty],\alpha\in\mathbb R$ and $\tau$ be an $\mathscr{F}_t$-stopping time. We say that ``equation (\ref{sto. NS}) driven by noise term (\ref{general noise form}) has a solution $u\in \widetilde{L}^p\big(0,\tau;\dot{B}_{q,\infty}^{\alpha}\big)$" if for every $\omega\in\Omega$, there exists $u^\omega\in \widetilde{L}^p\big(0,\tau(\omega);\dot{B}_{q,\infty}^{\alpha}\big)$ such that the equality
\begin{equation}
	\begin{aligned}
		u(t)= g(t)-\int_0^t\mathbf D_{s,t}(\varepsilon_{lin},{\rm B})e^{\nu(t-s)\Delta}\mathcal{P}{\rm div}\big(u(s)\otimes u(s)\big){\rm d}s\nonumber
	\end{aligned}
\end{equation}
holds in $\widetilde{L}^p\big(0,\tau;\dot{B}_{q,\infty}^{\alpha}\big)$ $\mathbb P$-almost surely. 
\end{definition}

\begin{remark}[About adaptiveness of solutions]
    Solutions constructed in the paper are based on standard fixed point argument built on continuity of the corresponding bilinear map and the classical contraction lemma (lemma \ref{Lem. contraction}). This means that one can always apply Picard iteration or Oseen iteration to get sequence that converges to the solution almost surely and strongly in the solution space. In this sense, the adaptiveness of the solution as an stochastic process can always be established. With this in mind, we would not be discussing the adaptiveness issue specially in this paper.
\end{remark}

\begin{remark}
	In deterministic case, i.e. when $\varepsilon_{lin}=0$ and $G_j\equiv0$, global well-posedness for small data are famously known in the space $C_{w^*,loc}\big([0,\infty);\dot{B}_{q,\infty}^{-1+d/q}\big)$ for $1\leq q<\infty$ \cite{CMP93,BCD11}. Such solutions are called ``mild solutions" as they preserve continuous trajectories in space $\dot{B}_{q,\infty}^{-1+d/q}$.
\end{remark}

Technically, in this paper, it is crucial to utilize time integral to achieve ``smallness". So we have to work on solution space with time integrability instead. For convenience, we denote the space for each $p\in[2,\infty]$ and $q\in [1,\infty]$ :
\begin{equation}
	\Xi_{p,q}(T):=\widetilde{L}^p\big(0,T;\dot{B}_{q,\infty}^{s+2/p}\big)\cap \widetilde{L}^{p^*}\big(0,T;\dot{B}_{q,\infty}^{s+2/p^*}\big)\nonumber
\end{equation}
with exponents
\begin{equation}\label{regularity exponents}
	p^*=\frac{p}{p-1}\quad and\quad s=s_q=-1+\frac{d}{q}~.
\end{equation}
In the space, we introduce the norm
\begin{equation}\label{norm of solution space}
	\|u\|_{\Xi_{p,q}(T)}:=\nu^{1/p}\|u\|_{\widetilde{L}^p(0,T;\dot{B}_{q,\infty}^{s+2/p})}+\nu^{1/p^*}\|u\|_{\widetilde{L}^{p^*}(0,T;\dot{B}_{q,\infty}^{s+2/p^*})}~.
\end{equation}
The following result is a slight generalization of the Cannone-Meyer-Planchon Theorem \cite{CMP93}. 

\begin{theorem}\label{Thm. existence}
	Let $\varepsilon_{lin}=0$ and $G_j\equiv0$, and $p\in(2,\infty]$, $q\in [1,\infty)$.
	 Given any initial value $u_0\in\dot{B}_{q,\infty}^{s}$, then at least for time scale $T^*\in(0,\infty]$ such that 
	\begin{equation}
		4C_1\nu^{-1}\big\|e^{\nu t\Delta}u_0\big\|_{\Xi_{p,q}(T^*)}\leq1~,\nonumber
	\end{equation}
	 the equation (\ref{sto. NS}) has a unique solution $u\in \Xi_{p,q}(T^*)$. The solution satisfies the bound:
	\begin{equation}
		\|u\|_{\Xi_{p,q}(T^*)}\leq 2~\big\|e^{\nu t\Delta}u_0\big\|_{\Xi_{p,q}(T^*)}~.\nonumber
	\end{equation}
	Here, $C_1=C_1(p,q,d)>0$ is the constant given by the continuity propostion \ref{Prop. continuity of bilinear map}. Moreover, 
	\begin{align}
		u\in \Xi_{\infty,q}(T^*)=L^\infty(0,T^*;\dot{B}_{q,\infty}^{-1+d/q}\big)\cap \widetilde{L}^{1}\big(0,T^*;\dot{B}_{q,\infty}^{1+d/q}\big)~.\nonumber
	\end{align}
\end{theorem}

\begin{corollary}
	There exists constant $C_0>0$ such that if
	\begin{equation}
		{\bf R}_q:=\nu^{-1}\|u_0\|_{\dot{B}_{q,\infty}^{s}}\leq C_0~,\nonumber
	\end{equation}
	the solution in Theorem \ref{Thm. existence} is global in time. For initial data with ${\bf supp}~\widehat{u}_0\subseteq \mathbb R^d\backslash B(0,r)$, the constrain can be relaxed slightly:
	\begin{equation}
		{\bf R}_q\leq C_0\min\big\{r^{2/p}, r^{2/p^*}\big\}.\nonumber
	\end{equation}
\end{corollary}

\begin{remark}
	The corollary simply follows from the inequality $\big\|e^{\nu t\Delta}u_0\big\|_{\Xi_{p,q}(\infty)}\lesssim \|u_0\|_{\dot{B}_{q,\infty}^{s}}$ for general $u_0\in\dot{B}_{q,\infty}^{s}$; and the inequality $\nu^{1/p}\big\|e^{\nu t\Delta}u_0\big\|_{\widetilde{L}^p(0,\infty ;\dot{B}_{q,\infty}^{s+2/p})}\lesssim r^{-2/p}\|u_0\|_{\dot{B}_{q,\infty}^{s}}$ for $u_0\in\dot{B}_{q,\infty}^{s}$ satisfying ${\bf supp}~\widehat{u}_0\subseteq \mathbb R^d\backslash B(0,r)$~. Both inequalities are straight-forward by using standard heat estimate for dyadic blocks (see Lemma \ref{Lem. heat estimates}).
\end{remark}

\
\subsubsection{Regularization by additive type noise ($\varepsilon_{lin}=0$)} It is straightforward to extend Theorem \ref{Thm. existence} onto stochastic setting for general additive noise. However, a crucial ``advantage" comparing to the deterministic side is that one could get a high probability of achieving unique global solution, at least if the noise is well-selected. Several choices of noise are come up with in this paper by utilizing the exponential martingale 
\begin{align}\label{exponential martingale}
	\mathbf M_t(\kappa,{\rm B}):=\exp\left(\kappa{\rm B}_t-\frac{\kappa^2}{2}t\right)\quad\quad 0\leq t<+\infty
\end{align}
where ${\rm B}$ is a one dimensional standard Brownian motion and $\kappa\neq 0$. We start with a simple case: if the initial value has compact support of Fourier frequency or more generally, belongs to $\dot{B}_{q,\infty}^{s}\cap\dot{B}_{q,\infty}^{s+}$, then the Stokes flow can be easily controlled by a single noise term:
\begin{align}\label{construction of noise I}
	\mathcal N(t)=\mathcal N^{add}(t)=\int_0^tG(s){\rm dB}_s~,\quad
	G(t)=G^{\kappa,N}(t)=\kappa\mathbf M_t(\kappa,{\rm B})\sum_{j\leq N}e^{\nu t\Delta}\dot{\Delta}_j u_0
\end{align}
where $\dot{\Delta}_j$ denote the homogenous dyadic blocks in Littlewood-Paley decomposition and $\kappa\neq0,N\in\mathbb N$ are constants to be determined. In the following we denote by $C=C(p,d)>0$ some constant depends only on $p$ and $d$.

\begin{theorem}\label{Thm. stochastic perturbation for slightly smoother initial data}
    Let $p\in(2,\infty)$, $q\in[1,\infty)$, $\theta\in(0,1)$ and $p^*,s$ be given by (\ref{regularity exponents}). Let also $u_0\in\dot{B}_{q,\infty}^{s}\cap\dot{B}_{q,\infty}^{s+\delta}$ with some $\delta>0$. Then for every
    \begin{equation}
		\begin{aligned}\label{constrain 1}
			N\geq \delta^{-1}\log_2\Big(C \nu^{-1}\|u_0\|_{\dot{B}_{q,\infty}^{s+\delta}}\Big) \quad and\quad \kappa^2>\frac{C\nu 2^{2N+1}}{1-\theta}\left({\bf R}_q^p+{\bf R}_q^{p^*}\right)=:K_q^2
		\end{aligned}
	\end{equation}
	there exists a Brownian stopping time $\tau$ such that: Equation (\ref{sto. NS}) driven by noise term (\ref{construction of noise I}) has a unique solution $u\in\Xi_{p,q}(\tau)$. 
   Moreover, the lower bound of probability holds:
	\begin{equation}
		\mathbb P\big(\tau=\infty\big)\geq 1-\left(\frac{\kappa}{K_q}\right)^{-\frac{2\theta}{p+p^*}}.\nonumber
	\end{equation} 
\end{theorem}

\

\noindent To treat general data in $\dot{B}_{q,\infty}^{s}$, we need extra noise terms to control the high frequency parts of the Stokes flow independently. In this case, they are given by
\begin{equation}\label{construction of noise II}
    \left\{\begin{aligned}
        \mathcal N(t)
	&=\mathcal N^{add}_{[{\rm B},{\rm B}^j]}(t)=\int_0^tG(s){\rm dB}_s+\sum_{j>N}\int_0^tG_j(s){\rm dB}_s^j\\	
	G(t)
	&=\kappa\mathbf M_t(\kappa,{\rm B})\sum_{j\leq N}e^{\nu t\Delta}\dot{\Delta}_j u_0\\
	G_j(t)
    &=\eta_j\mathbf M_t\left(\eta_j,{\rm B}^j\right)e^{\nu t\Delta}\dot{\Delta}_j u_0\quad (j>N)
    \end{aligned}\right.
\end{equation}
where $\left\{{\rm B},{\rm B}^j\right\}$ is an independent family of Brownian motions and
\begin{align}\label{eta_j}
		\eta_j=C\sqrt{{\bf R}_q^p+{\bf R}_q^{p^*}}~\nu^{1/2} 2^{j(1+\gamma)}~,\quad \gamma>0.
\end{align}

\begin{theorem}\label{Thm. stochastic perturbation II}
	Let the exponents be given as in Theorem \ref{Thm. stochastic perturbation for slightly smoother initial data} and $u_0\in \dot{B}_{q,\infty}^{s}$. Then for every large enough $N\in\mathbb N$ and $\kappa^2>K_q^2$, where $K_q>0$ is given as in (\ref{constrain 1}), there exists a Brownian stopping time $\tau$ such that: Equation (\ref{sto. NS}) driven by noise term (\ref{construction of noise II}) has a unique solution $u\in\Xi_{p,q}(\tau)$. 
    Moreover,
	\begin{align}
		\mathbb P\big(\tau=\infty\big)\geq\left[1-\left(\frac{\kappa}{K_q}\right)^{-\frac{2\theta}{p+p^*}}\right]\exp\Bigg(-\frac{6}{\pi^2}\sum_{j>N}j^{-2}\Bigg)~.\nonumber
	\end{align} 
\end{theorem}

\
\begin{remark}\label{remark. additional properties of the solution I}
    The solutions $u$ given by Theorem \ref{Thm. stochastic perturbation for slightly smoother initial data} and Theorem \ref{Thm. stochastic perturbation II} have continuous paths in the sense that $u\in C_{w^*,loc}\big([0,\tau);\dot{B}_{q,\infty}^{s}\big)$ and satisfies the bound $\|u\|_{\Xi_{p,q}(\tau)}\leq 2\|u_0\|_{\dot{B}_{q,\infty}^{s}}$~. Note that $g$ is given by (\ref{general noise form}) and, in general, solves the stochastic Stokes equations (\ref{sto. Stokes}) in the space $C_{loc}\big([0,\infty);\mathcal S'\big)$. And it is not hard to show that $g\in C_{loc}\big((0,\tau);\dot{B}_{q,\infty}^{s}\big)\cap C_{w^*,loc}\big([0,\tau);\dot{B}_{q,\infty}^{s}\big)$, given the stopping time defined as in subsection \ref{subsec. Domination by noise I}. However, strong $\dot{B}_{q,\infty}^{s}$-continuity near $t=0$ or even continuity in stronger topology would generally require either appropriate decaying of the noise coefficient as $|j|\rightarrow\infty$ (which is not the case in this paper) or further regularity condition on $u_0$.
\end{remark}

\
\subsubsection{Restoring well-posedness in $\dot{B}^{-1}_{\infty,\infty}$}
If we allow linear multiplicative noise in the noise term ( i.e. $\varepsilon_{lin}\neq0$ ), we would be able to restore a local well-posedness theory properly for initial values in the space $\dot{B}^{-1}_{\infty,\infty}$. The theory yields for every large enough $T>0$, that the solution lives on finite interval $[0,T]$ in high probability. In this part, we look at the solution space given by
\begin{align}
    X^{p,\delta}(T):=\Xi_{p,\infty}(T)\cap\widetilde{L}^{p}\big(0,T;\dot{B}_{\infty,\infty}^{-1+2/p+\delta}\big)\cap\widetilde{L}^{p^*}\big(0,T;\dot{B}_{\infty,\infty}^{-1+2/p^*+\delta}\big)~.\nonumber
\end{align} 
See (\ref{norm of space X}) for the norm of space $X^{p,\delta}(T)$. As before, the result depends on the the exact noise form. As the simplest case, we first investigate the case of a single linear multiplicative noise, i.e.
\begin{align}\label{construction of noise III}
    \mathcal N(t)=\varepsilon\nu^{1/2}\int_0^tu(s){\rm dB}_s~,\quad 0<\varepsilon^2<1~.
\end{align}
Note that in the general form (\ref{general noise form}), this corresponds to 
\begin{align}
    \mathcal N(t)=\varepsilon\nu^{1/2}\int_0^t\big(u(s)-g(s)\big){\rm dB}_s+\varepsilon\nu^{1/2}\int_0^tg(s){\rm dB}_s\nonumber
\end{align}
with $g(t)=\mathbf M_t(\varepsilon\nu^{1/2},{\rm B})~e^{\nu t\Delta}u_0$~. 

\begin{theorem}\label{Thm. restoring well-posedness by noise}
    Let $p\in(2,\infty)$ and $u_0\in \dot{B}_{\infty,\infty}^{-1}\cap\dot{B}_{\infty,\infty}^{-1+\delta}$ with $0<\delta<\min\{2/p,2/p^*\}$. Then for every $T>0$ there exists a Brownian stopping time $\mathcal T$ and constant $B=B(p,\delta,d)>0$ such that: For any initial value satisfying $\|u_0\|_{\dot{B}_{\infty,\infty}^{-1+\delta}}\leq B\nu\varepsilon^4e^{-\nu\varepsilon^2T/2}$, equation (\ref{sto. NS}) driven by noise term (\ref{construction of noise III}) has a unique solution $u\in X^{p,\delta}(\mathcal T)$. Moreover, there exist universal constants $r_0,c_0>0$ such that if $~T\geq r_0/(\nu\varepsilon^2)$, then
        \begin{align}
            \mathbb P\big(\mathcal T>T\big)\geq 1-2e^{-c_0\nu\varepsilon^2T}.\nonumber
        \end{align}
\end{theorem}

\begin{remark}
    We note here that by considering noise term given instead by
    \begin{align}
        \mathcal N(t)=\varepsilon\nu^{1/2}\int_0^t\big(u(s)-g(s)\big){\rm dB}_s+\kappa\int_0^tg(s){\rm dB}_s'\nonumber
    \end{align}
    with $g(t)=\mathbf M_t(\kappa,{\rm B})~e^{\nu t\Delta}\dot{S}_{N+1}u_0$ (${\rm B}$ and ${\rm B'}$ are independent), one can drop the smallness condition of $u_0$ in Theorem \ref{Thm. restoring well-posedness by noise}. This is achieved by making $N$ and $\kappa^2$ large enough in the spirit of Theorem \ref{Thm. stochastic perturbation for slightly smoother initial data}. The exact theorem formulation and the proof are left to the readers.
\end{remark}

\
\par Like in Theorem \ref{Thm. stochastic perturbation II}, we can remove the additional regularity assumption and the smallness condition of initial value at the same time by introducing independent additive noise terms that controls the associated Stokes flow. In the following theorem, our noise term $\mathcal N$ would be given by
\begin{equation}\label{construction of noise IV}
    \left\{\begin{aligned}
        &\mathcal N(t)
        =\varepsilon\nu^{1/2}\int_0^t\big(u(s)-g(s)\big){\rm dB}_s+\mathcal N^{add}(t)~,\quad 0<\varepsilon^2<1\\
        &g(t)
        =e^{\nu t\Delta}u_0+\int_0^te^{\nu(t-s)\Delta}{\rm d}\mathcal N^{add}(s)\\
        &\mathcal N^{add}(t)
        =\mathcal N^{add}_{[\widetilde{\rm B},{\rm B}^j]}(t)\quad via\quad (\ref{construction of noise II})
    \end{aligned}\right.
\end{equation}
where $\big\{{\rm B},\widetilde{\rm B},{\rm B}^j\big\}$ is an independent family of Brownian motions.

\begin{theorem}\label{Thm. restoring well-posedness by noise II}
    Let $p\in(2,\infty)$, $\theta\in(0,1)$, $0<\delta<\min\{2/p,2/p^*\}$ and $u_0\in \dot{B}_{\infty,\infty}^{-1}$. Then for every large enough $N\in\mathbb N_0$, every $T>0$ and
    \begin{align}
        \kappa^2
        &>\frac{C\nu 2^{4N+1}}{1-\theta}\left(\frac{e^{\nu\varepsilon^2T/4}}{\varepsilon^2}\right)^{p+p^*}\left({\bf R}_\infty^p+{\bf R}_\infty^{p^*}\right)=:K_\infty^2\nonumber\\
        \eta_j
        &=C\left(\frac{e^{\nu\varepsilon^2T/4}}{\varepsilon^2}\right)^{\frac{p+p^*}{2}}\sqrt{{\bf R}_\infty^p+{\bf R}_\infty^{p^*}}~\nu^{1/2} 2^{j(1+\gamma)}~,\quad \gamma>\frac{\delta}{2}(p+p^*)\nonumber
    \end{align}
    there exists a Brownian stopping time $\mathcal T$ such that: Equation (\ref{sto. NS}) driven by noise term (\ref{construction of noise IV}) has a unique solution $u\in X^{p,\delta}(\mathcal T)$. Moreover, there exist universal constants $r_0,c_0>0$ such that if $~T\geq r_0/(\nu\varepsilon^2)$, then
    \begin{align}
        \mathbb P\big(\mathcal T>T\big)\geq \left(1-e^{-c_0\nu\varepsilon^2T}\right)\left[1-\left(\frac{\kappa}{K_\infty}\right)^{-\frac{2\theta}{p+p^*}}\right]\exp\Bigg(-\frac{6}{\pi^2}\sum_{j>N}j^{-2}\Bigg)~.\nonumber
    \end{align}
\end{theorem}

\begin{remark}\label{remark. additional properties of the solution II}
    The solutions $u$ given by Theorem \ref{Thm. restoring well-posedness by noise} and Theorem \ref{Thm. restoring well-posedness by noise II} have continuous modified paths in the sense that $u=g+(u-g)\in C_{w^*,loc}\big([0,\mathcal T);\dot{B}_{\infty,\infty}^{-1}\big)+ C_{w^*,loc}\big([0,\mathcal T);\dot{B}_{\infty,\infty}^{-1+\delta}\big)$ and satisfies the bound $\|u\|_{X^{p,\delta}(\mathcal T)}\leq2\|g\|_{X^{p,\delta}(\mathcal T)}\lesssim\nu\varepsilon^2e^{-\nu\varepsilon^2T/4}$~.
\end{remark}

\
\subsection{Notation and some basic concepts}

Throughout the paper, we employ the notation $a\lesssim b$ if there exists a constant $C>0$ (independent of $a,b$) such that $a\leq Cb$. And by $C(p,q)$ we mean that the constant $C$ depends on parameters $p$ and $q$. For an integrability exponent $p\in[1,\infty]$, we always write $p^*:=\frac{p}{p-1}$ so that $\frac{1}{p}+\frac{1}{p^*}=1$. We would also use the notation $\mathbb N_0=\mathbb{N}\cup\{0\}$.\\

\noindent {\bf Balls in Banach space:} We denote by ${\rm B}\big(x,r;X\big)$ the open ball in Banach space $\big(X,\|\cdot\|_X\big)$ centered at $x\in X$ with radius $r>0$. For the closed ball, we simply write $\overline{{\rm B}\big(x,r;X\big)}$.\\

\noindent {\bf Dilation:} For (vector-valued) function $f:\mathbb R^d\rightarrow\mathbb R^n$, we denote its dilation $\delta_\lambda f(x):=f(\lambda x)$ with some $\lambda>0$.\\

\noindent {\bf Fourier Transforms:} In this paper, the Fourier transform and inverse Fourier transform of a function $f$ is defined by
\begin{equation}
	\mathcal F[f](\xi):=\int_{\mathbb R^d}f(x)e^{2\pi i\xi\cdot x}{\rm d}x,\quad
	\mathcal F^{-1}[f](\xi):=\int_{\mathbb R^d}f(x)e^{-2\pi i\xi\cdot x}{\rm d}x.\nonumber
\end{equation}
And we sometimes use the notations: $\widehat f=f^\wedge=\mathcal F[f]$, $\check f=f^\vee=\mathcal F^{-1}[f]$.\\

\noindent {\bf Sobolev Spaces:} Let $\mathcal{S}'$ be the space of tempered (or Schwartz) distributions on $\mathbb{R}^d$. For $s\in\mathbb R$ and $1\leq p\leq\infty$, we denote the Sobolev space $W^{s,p}=\big\{f\in\mathcal{S}':\big\|f\big\|_{W^{s,p}}:=\big\|(I-\Delta)^{s/2}f\big\|_{L^p}<\infty\big\}$. In case $p=2$, we use the more conventional notation $H^s=W^{s,2}$ and $\langle\cdot,\cdot\rangle_{H^s}$ for the corresponding inner product. Also, in case $s=0$, we have the usual Lebesgue spaces $L^p=W^{0,p}$. For simplicity, 
we also use the notation $L^2_\sigma=\left\{v\in L^2: {\rm div}v=0~in~weak~sense\right\}$.
\par We use a convenient notation for Sobolev interpolation: let $0<p<s<q$ with $s=\theta p+(1-\theta)q$ for some $\theta\in(0,1)$, we denote by ``$H^s=\big[H^p,H^q\big]_\theta$" the interpolation $\|f\|_{H^s}\leq\|f\|_{H^p}^{\theta}\|f\|_{H^q}^{1-\theta}$.\\

\noindent {\bf Besov Spaces:} Let $\mathcal{S}'_h$ be the space of distributions $f$ on $\mathbb{R}^d$ such that $\displaystyle\lim_{\lambda\rightarrow\infty}\|\theta(\lambda D)f\|_{L^\infty}=0$ for all $\theta\in C^\infty_c(\mathbb R^d)$. Note that we have
\begin{equation}
	S'_h=\left\{f\in S':~\exists~\theta\in C^\infty_c~s.t.~\theta(0)\neq0~and~~\lim_{\lambda\rightarrow\infty}\|\theta(\lambda D)f\|_{L^\infty}=0\right\}.\nonumber
\end{equation}
We would work on the homogenous Besov spaces defined by means of the following homogenous dyadic partition of unity.
\begin{lemma}[Dyadic Partition of Unity, \cite{BCD11}]
	Let $\mathcal C=\left\{\xi\in\mathbb R^d:~3/4\leq|\xi|\leq8/3\right\}$. There exists a radial function $\phi\in C^\infty_c(\mathcal C;[0,1])$ satisfying the following:
	\begin{align}
		&i).~\sum_{j\in\mathbb Z}\phi(2^{-j}\xi)=1,~\forall\xi\in \mathbb R^d\backslash\{0\};\nonumber\\
		&ii).~1/2\leq\sum_{j\in\mathbb Z}\phi(2^{-j}\xi)^2\leq1,~\forall\xi\in \mathbb R^d\backslash\{0\};\nonumber\\
		&iii).~{\bf supp}(\delta_{2^{-j}}\phi)\cap{\bf supp}(\delta_{2^{-j'}}\phi)\neq\emptyset~only~for~|j-j'|\leq1~.\label{intersect of supports}	
	\end{align}
\end{lemma}
\noindent The so-called homogenous dyadic blocks $\dot\Delta_j$ are operators defined by the associated Fourier multipliers:
\begin{align}
	\dot\Delta_j f:=\phi(2^{-j}D)f=(\delta_{2^{-j}}\phi\cdot \widehat{f}~)^\vee,\quad f\in S'.\nonumber
\end{align}
Then the Littlewood-Paley decomposition holds for $f\in S'_h$: $\displaystyle f=\sum_{j\in\mathbb Z}\dot\Delta_j f$. For $s\in\mathbb R$ and $q,r\in[1,\infty]$, the homogenous Besov space $\dot{B}_{q,r}^{s}(\mathbb R^d)$ is defined by
\begin{equation}
	\dot{B}_{q,r}^{s}=\Bigg\{f\in S'_h:~\|f\|_{\dot{B}_{q,r}^{s}}:=\Bigg(\sum_{j\in\mathbb Z}\big(2^{js}\|\dot\Delta_j f\|_{L^q}\big)^r\Bigg)^{1/r}<\infty\Bigg\}.\nonumber
\end{equation}
\noindent Note that when ``$s<d/q$" or ``$s=d/p,~r=1$" $(s,q,r)$ satifies
\begin{equation}\label{scales condition for completion of Besov space}
	\left\{\begin{aligned}
		&s<d/q~,\quad 1<r\leq\infty\\
		&s=d/q~,\quad r=1
	\end{aligned}\right.
\end{equation}
the normed space $\dot{B}_{q,r}^{s}$ is a Banach space while in general it is not the case \cite{BCD11}. \\

\noindent {\bf Function Spaces involving Time:} Given a Banach space $\big(Y,~\|\cdot\|_Y\big)$ and $I\subset \mathbb{R}$, we write $L^p(I;Y)$ for the space of $L^p$-integrable functions from $I$ to $Y$, equipped with the usual $L^p$-norm. When $I=(a,b)$ we simply write $L^p(a,b;Y)$.
We also use $L^p_{loc}(0,\infty;Y)$ to denote the space of functions $f$ from $[0,\infty)$ to $Y$ satisfying $f\in L^p(0,T;Y)$ for all $T>0$. Similarly, we write $C(I;Y)$ for the space of continuous functions from $I$ to $Y$ equipped with the supremum norm in a bounded subset. 
We denote by 
$C_{w^*}(I;Y^*)$ the space of functions $u: I\rightarrow Y^*$ such that for each $v\in Y$, $t\mapsto\langle u(t),v\rangle$ is a real continuous function. We would also use $C_{w^*,loc}(I;Y^*)$ to denote functions $u: I\rightarrow Y^*$ such that $u\in C_{w^*}(J;Y^*)$ for every compact $J\subset I$.
\par As in Besov estimates it is usually easier to estimate each dyadic block in $L^p(0,T;L^q)$, we will find it useful to define the norm
\begin{equation}
	\|u\|_{\widetilde L^p(0,T;\dot{B}_{q,r}^{s})}:=\Bigg(\sum_{j\in\mathbb Z}\left(2^{js}\|\dot\Delta_j u\|_{L^p(0,T;L^q)}\right)^r\Bigg)^{1/r}\nonumber
\end{equation}
and the space $\widetilde L^p\big(0,T;\dot{B}_{q,r}^{s}\big)$ of tempered distributions $u$ on $(0,T)\times\mathbb R^d$ such that $\displaystyle\lim_{n\rightarrow-\infty}\big\|\sum_{j\leq n}\dot\Delta_j u\big\|_{L^p(0,T;L^\infty)}\\=0$ and $\displaystyle\|u\|_{\widetilde L^p(0,T;\dot{B}_{q,r}^{s})}<\infty$. It's not hard to check that when $(s,q,r)$ satisfies (\ref{scales condition for completion of Besov space}), the space $\widetilde L^p\big(0,T;\dot{B}_{q,r}^{s}\big)$ is a Banach space. From Minkowski inequality, one sees that
\begin{equation}
	\left\{\begin{aligned}
		&\|u\|_{\widetilde L^p(0,T;\dot{B}_{q,r}^{s})}\leq\|u\|_{L^p(0,T;\dot{B}_{q,r}^{s})}\quad when~p\leq r~;\\
		&\|u\|_{L^p(0,T;\dot{B}_{q,r}^{s})}\leq\|u\|_{\widetilde L^p(0,T;\dot{B}_{q,r}^{s})}\quad when~r\leq p~.\nonumber
	\end{aligned}\right.
\end{equation}

\subsection{Organization}
I will return to this part.

\
\section{Solutions in space $\Xi_{p,q}(T)$ for $1\leq q<\infty$}

\par This section is devoted to Theorem \ref{Thm. existence}, \ref{Thm. stochastic perturbation for slightly smoother initial data} and \ref{Thm. stochastic perturbation II} which are based on the following continuity proposition for the bilinear map 
\begin{equation}\label{bilinear mapping}
	\mathcal B_t(u,v):=-\int_0^te^{\nu(t-s)\Delta}\mathcal{P}{\rm div}\big(u(s)\otimes v(s)\big){\rm d}s~,\quad 0\leq t<\infty.
\end{equation}
As in the two theorems our noise terms do not rely on the solution, we can write the equation (\ref{sto. NS}) in the abstract form:
\begin{align}\label{sto. NS abstract I}
    u=g+\mathcal B(u,u)
\end{align}
where $g=e^{\nu t\Delta}u_0+\int_0^te^{\nu (t-s)\Delta}{\rm d}\mathcal N(s)$~.

\begin{proposition}\label{Prop. continuity of bilinear map}
	Let $r\in[2,\infty]$, $p\in(2,\infty]$ and $q\in[1,\infty)$. Then for any $0<T\leq\infty$, the bilinear mapping $\mathcal B: \Xi_{p,q}(T)\times \Xi_{p,q}(T)\rightarrow \Xi_{r,q}(T)$ is continuous and we have
    \begin{align}
        \|\mathcal B(u,v)\|_{\Xi_{r,q}(T)}\leq C_1 \nu^{-1}\|u\|_{\Xi_{p,q}(T)}\|v\|_{\Xi_{p,q}(T)}~,\quad u,v\in \Xi_{p,q}(T)~.\nonumber
    \end{align}
    In particular,
	\begin{align}
		\|\mathcal B(u,v)\|_{\Xi_{p,q}(T)}
        &\leq C_1 \nu^{-1}\|u\|_{\Xi_{p,q}(T)}\|v\|_{\Xi_{p,q}(T)}~,\label{continuity of bilinear map}\\
        \|\mathcal B(u,v)\|_{\Xi_{\infty,q}(T)}
        &\leq C_1 \nu^{-1}\|u\|_{\Xi_{p,q}(T)}\|v\|_{\Xi_{p,q}(T)}~.\label{bilinear estimate for mild space}
	\end{align}
    Furthermore, for each $\varphi\in\dot{B}_{q^*,1}^{1-d/q}$, the function $t\mapsto \big\langle B_t(u,v),\varphi\big\rangle$ is continuous on $[0,T]$ provided that $u,v\in \Xi_{p,q}(T)$.
\end{proposition}

\noindent As well as we know, Proposition \ref{Prop. continuity of bilinear map} is not included in any of the literature. However, its proof is different from its classical version (continuity on space $\Xi_{\infty,q}(T)$) just in the sense that one needs to apply Young's inequality more skillfully and choose the right regularity exponents to make sums converge. We leave the proof in subsection \ref{sec. Continuity of the Bilinear map} for convenience. As in the classical theory of well-posedness with smallness condition, one can derive ``small" solutions in the space $\Xi_{p,q}(T)$ by Proposition \ref{Prop. continuity of bilinear map} and the following contraction lemma, given that \begin{align}\label{smallness of zeroth approximation}
    4C_1\nu^{-1}\|g\|_{\Xi_{p,q}(T)}<1.
\end{align} 
Here $C_1$ is just the constant in (\ref{continuity of bilinear map}).

\begin{lemma}[see \cite{CG09}]\label{Lem. contraction}
	Let $E$ be a Banach space and $\mathcal B$ a continuous bilinear map on $E$. Then for any $g\in\overline{{\rm B}(0,R;E)}$ with $0<R<\frac{1}{4\|\mathcal B\|}$, there exists a unique $x\in \overline{{\rm B}(0,2R;E)}$ such that $x=g+\mathcal B(x,x)$. 
\end{lemma}

\begin{proposition}\label{Prop. small solutions to the abstract equation}
    The abstract equation (\ref{sto. NS abstract I}) has a unique solution $u\in\overline{{\rm B}\big(0,2\|g\|;~\Xi_{p,q}(T)\big)}$ when the smallness condition (\ref{smallness of zeroth approximation}) holds. 
\end{proposition}
\begin{proof}[Proof of Proposition \ref{Prop. small solutions to the abstract equation}]
    It follows directly by Proposition \ref{Prop. continuity of bilinear map} and lemma \ref{Lem. contraction}.
\end{proof}
\noindent In the next subsection, we proof Theorem \ref{Thm. stochastic perturbation for slightly smoother initial data} and \ref{Thm. stochastic perturbation II} respectively by constructing appropriate stopping time $\tau$ up to which (\ref{smallness of zeroth approximation}) holds. To achieve this, it is essential to work on the spaces $\Xi_{p,q}(\tau)$ ($p\neq\infty$) rather than the classical $\Xi_{\infty,q}(T)$. After constructing solution $u\in \Xi_{p,q}(\tau)$, (\ref{bilinear estimate for mild space}) then implies that the solution paths $u=g+(u-g)$ also lay in the mild path space $\Xi_{\infty,q}=L^\infty_t\dot{B}_{q,\infty}^{-1+d/q}\cap\widetilde L^1_t\dot{B}_{q,\infty}^{1+d/q}$. Also, it is not hard to see that $g$, $u-g$ and hence $u$ belongs to $C_{w^*,loc}\big([0,\tau);\dot{B}_{q,\infty}^{-1+d/q}\big)$.

\
\subsection{Domination by noise}\label{subsec. Domination by noise I}
Note that with noise term $\mathcal N$ given by (\ref{construction of noise I}), one can write
\begin{equation}
	\begin{aligned}
	g(t)=\underbrace{\sum_{j\leq N}\mathbf M_t(\kappa,{\rm B})~e^{\nu t\Delta}\dot\Delta_j u_0}_{\mathcal L^N}+\underbrace{\sum_{j> N}e^{\nu t\Delta}\dot\Delta_j u_0}_{\mathcal H^N}~;\label{split the zeroth approximation}
\end{aligned}
\end{equation}
and if $\mathcal N$ is given by (\ref{construction of noise II}), one writes instead:
\begin{equation}
	\begin{aligned}
	g(t)=\sum_{j\leq N}\mathbf M_t(\kappa,{\rm B})~e^{\nu t\Delta}\dot\Delta_j u_0+\underbrace{\sum_{j> N}\mathbf M_t\left(\eta_j,{\rm B}^j\right)e^{\nu t\Delta}\dot\Delta_j u_0}_{\widetilde{\mathcal H}^N}~.\label{split the zeroth approximation II}
\end{aligned}
\end{equation}
To control the exponential martingale $\mathbf M_t(\kappa,{\rm B})$, we introduce the stopping time for $0<\theta<1$ and $A>1$:
\begin{align}\label{stopping time to control M_t}
    \tau^{\theta,\kappa}_A({\rm B}):=\inf\Bigg\{t>0~:~\exp\left(\kappa{\rm B}_t-\theta\frac{\kappa^2}{2}t\right)>A\Bigg\}~.
\end{align}

\noindent We will show the following two propositions which indicate Theorem \ref{Thm. stochastic perturbation for slightly smoother initial data} and \ref{Thm. stochastic perturbation II} respectively.

\begin{proposition}\label{Prop. smallness of g I}
    Let $p\in[2,\infty),q\in[1,\infty],\theta\in(0,1)$ and $s=-1+d/q$. Let also $g$ be given by (\ref{split the zeroth approximation}) with $u_0\in\dot{B}_{q,\infty}^{s}\cap\dot{B}_{q,\infty}^{s+\delta}$ and $N,\kappa$ satisfying (\ref{constrain 1}). Take $A=(\kappa/K)^{\frac{2}{p+p^*}}$, then up to $\tau=\tau^{\theta,\kappa}_A({\rm B})$, we have that
    $4C_1\nu^{-1}\|g\|_{\Xi_{p,q}(\tau)}<1$. Moreover,
    \begin{equation}\label{lower bound of probability I}
		\mathbb P\big(\tau=\infty\big)\geq 1-\left(\frac{\kappa}{K}\right)^{-\frac{2\theta}{p+p^*}}.
	\end{equation} 
\end{proposition}

\begin{proposition}\label{Prop. smallness of g II}
    Let the exponents $p,q,s,\theta$ be given as in Proposition \ref{Prop. smallness of g I} and $u_0\in\dot{B}_{q,\infty}^{s}$. Let also $g$ be given by (\ref{split the zeroth approximation II}) with $\kappa$ satisfying (\ref{constrain 1}); $\eta_j$ given by (\ref{eta_j}) and $N\geq N_0$ for some constant $N_0>0$. Take $A=(\kappa/K)^{\frac{2}{p+p^*}}$, $A_j=\left(\frac{\pi^2}{12}j^2\right)^2$ and
    \begin{align}
        \tau=\tau^{\theta,\kappa}_A({\rm B})\wedge\bigwedge_{j>N}\tau^{1/2,\eta_j}_{A_j}({\rm B}^j)~.\nonumber
    \end{align}
    Then we have $4C_1\nu^{-1}\|g\|_{\Xi_{p,q}(\tau)}<1$ and
    \begin{align}
		\mathbb P\big(\tau=\infty\big)\geq\left[1-\left(\frac{\kappa}{K}\right)^{-\frac{2\theta}{p+p^*}}\right]\exp\Bigg(-\frac{6}{\pi^2}\sum_{j>N}j^{-2}\Bigg)~.\nonumber
	\end{align} 
\end{proposition}

\noindent The proofs are based on the following lemma whose proof we delay until the end of this subsection.

\
\begin{lemma}\label{Lem. probability I}
	Let $\tau^{\theta,\kappa}_A({\rm B})$ be given by (\ref{stopping time to control M_t}). Then we have $\mathbb P\left(\tau^{\theta,\kappa}_A({\rm B})=\infty\right)\geq 1-A^{-\theta}$.
\end{lemma}

\
\begin{proof}[Proof of Proposition \ref{Prop. smallness of g I}]
    In what follows, we denote by $C=C(p,d)>0$ some constant that could be different from line to line. We remind the splitting $g=\mathcal L^N+\mathcal H^N$ given by (\ref{split the zeroth approximation}). The control of the norm would be given respectively.

	\noindent {$\bullet~\bf High~Frequency~part:$} By standard heat estimate for dyadic blocks given in Lemma \ref{Lem. heat estimates}, we have that
	\begin{align}
		\big\|e^{\nu t\Delta}\dot\Delta_j u_0\big\|_{L^p(0,\infty; L^q)}\lesssim \nu^{-\frac{1}{p}}2^{-j\frac{2}{p}}\big\|\dot{\Delta}_j u_0\big\|_{L^q}\nonumber
	\end{align}
	and then
	\begin{align}
		\nu^{\frac{1}{p}}2^{j\left(s+\frac{2}{p}\right)}\big\|e^{\nu t\Delta}\dot\Delta_j u_0\big\|_{L^p(0,\infty; L^q)}\lesssim 2^{-j\delta}\|u_0\|_{\dot{B}_{q,\infty}^{s+\delta}}\nonumber
	\end{align}
	with $s=-1+\frac{d}{q}$. This gives that
	\begin{align}
		\nu^{\frac{1}{p}}\big\|\Delta_j\mathcal H_N\big\|_{L^p(0,\infty; L^q)}
		\lesssim \sum_{\ell=j-1}^{j+1}2^{-\ell\left(s+\frac{2}{p}+\delta\right)}\|u_0\|_{\dot{B}_{q,\infty}^{s+\delta}}\lesssim 2^{-j\left(s+\frac{2}{p}+\delta\right)}\|u_0\|_{\dot{B}_{q,\infty}^{s+\delta}}\nonumber
	\end{align}
	and hence
	\begin{align}
		\nu^{\frac{1}{p}}\big\|\mathcal H_N\big\|_{\widetilde{L}^p\big(0,\infty; \dot{B}_{q,\infty}^{s+2/p}\big)}
		\lesssim 2^{-N\delta}\|u_0\|_{\dot{B}_{q,\infty}^{s+\delta}}~.\nonumber
	\end{align}
	And by definition (\ref{norm of solution space}) of $\Xi_{p,q}$-norm,
	\begin{align}\label{control of high freq. via loss of regularity}
		4\nu^{-1}\big\|\mathcal H_N\big\|_{\Xi_{p,q}(\infty)}\leq C 2^{-N\delta}\nu^{-1}\|u_0\|_{\dot{B}_{q,\infty}^{s+\delta}}\leq 1/4~,
	\end{align}
	as long as we choose a sufficiently large $N\in\mathbb N$ such that
	\begin{align}
		N\geq \delta^{-1}\log_2\left(\frac{C\|u_0\|_{\dot{B}_{q,\infty}^{s+\delta}}}{\nu}\right)~.\nonumber
	\end{align}
	
	\noindent {$\bullet~\bf Low~Frequency~part:$} Let $\theta\in(0,1)$ be arbitrarily given, $\kappa\neq0$ and $A>1$ be determined later. Then up to the stopping time $\tau=\tau^{\theta,\kappa}_A$, we obviously have $\mathbf M_t(\kappa,{\rm B})\leq Ae^{-(1-\theta)\frac{\kappa^2}{2}t}$. Together by lemma \ref{Lem. heat estimates}, we derive that
    \begin{align}
        \big\|\mathbf M_t(\kappa,{\rm B})e^{\nu t\Delta}\dot\Delta_j u_0\big\|_{L^q}\lesssim Ae^{-t\left((1-\theta)\frac{\kappa^2}{2}+c\nu 2^{2j}\right)}\|\dot\Delta_j u_0\big\|_{L^q}\nonumber
    \end{align}
    and so
	\begin{align}\label{low freq. blocks estimate of g}
		4\nu^{\frac{1}{p}-1}2^{j\left(s+\frac{2}{p}\right)}\big\|\dot\Delta_j \mathcal L^N\big\|_{L^p(0,\tau; L^q)}< C\left[\frac{\nu 2^{2j}}{(1-\theta)\frac{\kappa^2}{2}}\right]^{\frac{1}{p}}A\nu^{-1}\|u_0\|_{\dot{B}_{q,\infty}^{s}}
	\end{align}
	holds for all $j\leq N+1$. Now, if $\kappa\neq0$ and $A>1$ satisfy
	\begin{align}
		(1-\theta)\frac{\kappa^2}{2}\geq\left(\frac{4CA\|u_0\|_{\dot{B}_{q,\infty}^{s}}}{\nu}\right)^p \nu 2^{2N}\nonumber
	\end{align}
	so that $\displaystyle C\left[\frac{\nu 2^{2j}}{(1-\theta)\frac{\kappa^2}{2}}\right]^{\frac{1}{p}}A\nu^{-1}\|u_0\|_{\dot{B}_{q,\infty}^{s}}\leq 1/4$ for all $j\leq N+1$, we would see that
	\begin{align}
		4\nu^{-1}\big\|\mathcal L^N\big\|_{\widetilde{L}^p\big(0,\tau; \dot{B}_{q,\infty}^{s+2/p}\big)}\leq 1/4~.\nonumber
	\end{align}
	The same thing also holds with $p$ replaced by $p^*$, so if we set $\kappa\neq0$ and $A>1$ such that
	\begin{align}\label{lower bound condition for noise intensity}
		(1-\theta)\frac{\kappa^2}{2}\geq A^{p+p^*}\left[\left(\frac{4C\|u_0\|_{\dot{B}_{q,\infty}^{s}}}{\nu}\right)^p+\left(\frac{4C\|u_0\|_{\dot{B}_{q,\infty}^{s}}}{\nu}\right)^{p^*}\right] \nu 2^{2N}~,
	\end{align}
	we would achieve the control:
	\begin{align}\label{control of low frequency part}
		4\nu^{-1}\big\|\mathcal L^N\big\|_{\Xi_{p,q}(\tau)}\leq 1/4+1/4=1/2~.
	\end{align}
	Together with (\ref{control of high freq. via loss of regularity}), we have that $4C_1\nu^{-1}\|g\|_{\Xi_{p,q}(\tau)}<1$. To fulfill (\ref{lower bound condition for noise intensity}), we may just set
    \begin{align}
        A=\left(\frac{\kappa^2}{K_q^2}\right)^\frac{1}{p+p^*}\nonumber
    \end{align}
    for any pre-given $\kappa^2>K_q^2$ where
	\begin{align}
		K_q^2:=\frac{C\nu 2^{2N+1}}{1-\theta}\left({\bf R}_q^p+{\bf R}_q^{p^*}\right),\quad\quad 
		\left({\bf R}_q:=\frac{\|u_0\|_{\dot{B}_{q,\infty}^{s}}}{\nu}\right)\nonumber
	\end{align}
	Finally, (\ref{lower bound of probability I}) follows by lemma \ref{Lem. probability I}.
\end{proof}

\
\begin{proof}[Proof of Proposition \ref{Prop. smallness of g II}]
	In the setting of Proposition \ref{Prop. smallness of g II}, $g$ is given by (\ref{split the zeroth approximation II}), and the stopping time
    \begin{align}
        \tau=\tau^{\theta,\kappa}_A({\rm B})\wedge\bigwedge_{j>N}\tau^{1/2,\eta_j}_{A_j}({\rm B}^j)~.\nonumber
    \end{align}
    Control of the low frequency part can be given as in the previous proof by (\ref{control of low frequency part}) for any splitting number $N\in\mathbb N$. And so we again have $\kappa^2>K_q^2$ and $A=(\kappa/K_q)^{\frac{2}{p+p^*}}$ where $K_q^2$ is given as before. Now it remains to set appropriate $\eta_j$ and $A_j$, and bound the high frequency part. Remind that the constant $C=C(p,d)>0$ could be different from line to line.
    \par We first set the sequence $A_j:=\left(\frac{\pi^2}{12}j^2\right)^2$ so that $\displaystyle\sum_{j=1}^\infty A_j^{-1/2}=2$. Similar to (\ref{low freq. blocks estimate of g}), by lemma \ref{Lem. heat estimates} and the fact that $\mathbf M_t(\eta_j,{\rm B}^j)\leq A_je^{-\frac{\eta_j^2}{4}t}$ up to stopping time $\tau$, we have:  
	\begin{align}\label{control of high frequency part}
		4\nu^{\frac{1}{p}-1}2^{j\left(s+\frac{2}{p}\right)}\big\|\mathbf M\left(\eta_j,{\rm B}^j\right)e^{\nu t\Delta}\dot\Delta_j u_0\big\|_{L^p(0,\tau; L^q)}\leq C\left(\frac{\nu 2^{2j}}{\eta_j^2}\right)^{\frac{1}{p}}A_j\nu^{-1}\|u_0\|_{\dot{B}_{q,\infty}^{s}}
	\end{align}
	for each $j\geq N$. 
    Of course, we should at least set $N\in\mathbb N$ such that $A_N>1$. Now if we take non-decreasing $\eta_j$ such that
	\begin{align}\label{eta_j condition}
		\eta_{j-1}^2\geq (8C{\bf R}_q)^p\nu 2^{2j}\left(\frac{\pi^2}{12}(j+1)^2\right)^{2p},
	\end{align}
	(\ref{control of high frequency part}) would give
	\begin{align}
		4\nu^{\frac{1}{p}-1}2^{j\left(s+\frac{2}{p}\right)}\big\|\dot\Delta_j \widetilde{\mathcal H}^N\big\|_{L^p(0,\tau; L^q)}
		&\lesssim\sum_{\ell=j-1}^{j+1}\left(\frac{\nu 2^{2\ell}}{\eta_\ell^2}\right)^{\frac{1}{p}}A_\ell\nu^{-1}\|u_0\|_{\dot{B}_{q,\infty}^{s}}\nonumber\\
		&\leq C\left(\frac{\nu 2^{2j}}{\eta_{j-1}^2}\right)^{\frac{1}{p}}A_{j+1}\nu^{-1}\|u_0\|_{\dot{B}_{q,\infty}^{s}}\leq 1/8\nonumber
	\end{align}
	and hence
	\begin{align}
		4\nu^{-1}\big\|\widetilde{\mathcal H}^N\big\|_{\widetilde{L}^p\big(0,\tau; \dot{B}_{q,\infty}^{s+2/p}\big)}\leq 1/8~.\nonumber
	\end{align}
	The same thing is true when we replace $p$ by $p^*$. So in view of (\ref{eta_j condition}), we may set 
	\begin{align}
		\eta_j:=C\sqrt{{\bf R}_q^p+{\bf R}_q^{p^*}}\nu^{1/2} 2^{j(1+\gamma)}~,\quad \gamma>0\nonumber
	\end{align}
	and $N\in\mathbb N$ large enough so that 
    \begin{align}
        2^{\gamma j}\geq\left(\frac{\pi^2}{12}(j+1)^2\right)^{p+p^*}~,\quad \forall~j>N.\nonumber
    \end{align}
    Then we see that $4\nu^{-1}\big\|\widetilde{\mathcal H}^N\big\|_{\Xi_{p,q}(\tau)}\leq 1/8+1/8=1/4$. Together with the low frequency part (\ref{control of low frequency part}), we derive
	\begin{align}
		4\nu^{-1}\|g\|_{\Xi_{p,q}(\tau)}\leq 1/2+1/4=3/4<1.\nonumber
	\end{align}
	on the final subset $\displaystyle\bigcap_{j>N}\Omega_j$. To this end, we denote
    \begin{align}
        \Omega^{\theta,\kappa}
        =\Big\{\tau^{\theta,\kappa}_A({\rm B})=\infty\Big\}\quad and\quad
        \Omega_j
        =\Big\{\tau^{1/2,\eta_j}_{A_j}({\rm B}^j)=\infty\Big\}~,\quad j>N~.\nonumber
    \end{align}
    As in the previous proof, we now have
    \begin{align}
        \mathbb P\big(\Omega^{\theta,\kappa}\big)\geq 1-\left(\frac{\kappa}{K_q}\right)^{-\frac{2\theta}{p+p^*}}~.\nonumber
    \end{align}
    Next, lemma \ref{Lem. probability I} and independence of the Brownian motions give that
	\begin{align}
		\mathbb P\Bigg(\bigcap_{j>N}\Omega_j\Bigg)
		&\geq\prod_{j>N}\Big(1-A_j^{-1/2}\Big)=\exp\left[\sum_{j>N}\ln\Big(1-A_j^{-1/2}\Big)\right]\nonumber\\
		&\geq\exp\Bigg(-2\sum_{j>N}A_j^{-1/2}\Bigg)
		=\exp\Bigg(-\frac{6}{\pi^2}\sum_{j>N}j^{-2}\Bigg)~.\nonumber
	\end{align}
	Here we have used the elementary inequality $\ln(1-x)\geq -2x$ for $0\leq x\leq 1/2$ and we should also set $N$ large enough so that $A_N\geq 4$. Putting together, independence implies that 
	\begin{align}
	    \mathbb P\Bigg(\Omega^{\theta,\kappa}\cap\bigcap_{j>N}\Omega_j\Bigg)\geq\left[1-\left(\frac{\kappa}{K_q}\right)^{-\frac{2\theta}{p+p^*}}\right]\exp\Bigg(-\frac{6}{\pi^2}\sum_{j>N}j^{-2}\Bigg)~.\nonumber
	\end{align}
	This ends the proof.
\end{proof}

\
\begin{proof}[Proof of Lemma \ref{Lem. probability I}]
	The lemma is proved in \cite{GHV14} in different parameter setting, we write down the proof here for the sake of inclusion. Denote
	\begin{align}
		\mathcal M_{\theta,\kappa}(t)=\exp\left(\kappa{\rm B}_t-\theta\frac{\kappa^2}{2}t\right)~,\quad\quad 0\leq t<\infty~.\nonumber
	\end{align}
	Then clearly $\mathbb P\big(\Omega^{\theta,\kappa}_A\big)=\mathbb P\big(\tau_A=\infty\big)$ where $\tau_A$ is the stopping time given by
	\begin{align}
		\tau_A=\tau_A^{\theta,\kappa}:=\inf\Big\{t>0~:~\mathcal M_{\theta,\kappa}(t)>A\Big\}~.\nonumber
	\end{align}
	By It\^o's formula, we have for any $\lambda>0$ that
	\begin{align}
		\mathcal M_{\theta,\kappa}(t)^\lambda=1+\lambda\kappa\int_0^t\mathcal M_{\theta,\kappa}(s)^{\lambda}{\rm dB}_s+(\lambda-\theta)\lambda\frac{\kappa^2}{2}\int_0^t\mathcal M_{\theta,\kappa}(s)^{\lambda}{\rm d}s~.\nonumber
	\end{align}
	Then in particular, when $\lambda=\theta$, one sees that
	\begin{align}
		\mathbb E\Big[\mathcal M_{\theta,\kappa}\left(t\wedge\tau_A\right)^\theta\Big]=1~.\nonumber
	\end{align}
	Thus, by monotonicity of the sets $\big\{\tau_A\geq t\big\}$ and Markov inequality,
	\begin{align}
		\mathbb P\big(\tau_A=\infty\big)
		&=\lim_{t\rightarrow\infty}\mathbb P\big(\tau_A\geq t\big)\nonumber\\
		&=\lim_{t\rightarrow\infty}\Big[1-\mathbb P\Big(\mathcal M_{\theta,\kappa}\left(t\wedge\tau_A\right)>A\Big)\Big]\nonumber\\
		&\geq 1-A^{-\theta}~.\nonumber
	\end{align}
\end{proof}

\
\subsection{Continuity of the Bilinear map}\label{sec. Continuity of the Bilinear map}

\begin{proof}[Proof of Proposition \ref{Prop. continuity of bilinear map}]
	Let $u,v\in\dot{B}_{q,\infty}^\alpha$ or just $\mathcal S$ and remind that
	\begin{equation}
		\mathcal B_t(u,v):=-\int_0^te^{\nu(t-s)\Delta}\mathcal{P}{\rm div}\big(u(s)\otimes v(s)\big){\rm d}s\nonumber.
	\end{equation}
	We recall the notation $\displaystyle\dot{S}_nf:=\sum_{j\leq n-1}\dot{\Delta}_jf$ and write the decomposition
	\begin{align}\label{Bony decomposition}
		\mathcal B(u,v)
		&=\sum_{j\in\mathbb Z}\mathcal B\left(\dot{S}_{j} u, \dot{\Delta}_j v\right)+\sum_{j\in\mathbb Z}\mathcal B \left(\dot{\Delta}_j u, \dot{S}_{j+1} v\right) \nonumber\\
		&=\underbrace{\sum_{\ell\in\mathbb Z}\sum_{j\geq\ell-N_0}\dot{\Delta}_\ell~\mathcal B\left(\dot{S}_{j} u, \dot{\Delta}_j v\right)}_{=:\mathcal Q^1(u,v)}
		+\underbrace{\sum_{\ell\in\mathbb Z}\sum_{j\geq\ell-N_0}\dot{\Delta}_\ell~\mathcal B \left(\dot{\Delta}_j u, \dot{S}_{j+1} v\right)}_{=:\mathcal Q^2(u,v)}~.
	\end{align}
	As the terms $\mathcal Q^1$ and $\mathcal Q^2$ can be dealt  with in very similar way, we only prove the conclusion for $\mathcal Q^1$ here.
	
	\noindent$\bf{\bullet~Estimate~of~the~norm:~}$ Let $r,p\in[2,\infty]$ and $q\in[1,\infty)$ for now. By Lemma \ref{Lem. heat estimates} and \ref{Lem. homogeneous multipliers}, standard Bernstein inequalities and H\"older inequality,
	\begin{align}
		\left\|\dot{\Delta}_\ell~\mathcal B_t\left(\dot{S}_{j} u, \dot{\Delta}_j v\right)\right\|_{L^q}\lesssim \int_0^te^{-c\nu 2^{2\ell}(t-s)}~2^\ell~\big\|\dot{S}_{j} u(s)\big\|_{L^\infty}\big\|\dot{\Delta}_{j} v(s)\big\|_{L^q}{\rm d}s~.\nonumber
	\end{align}
	As Bernstein inequalities give $\displaystyle\big\|\dot{S}_{j} u(s)\big\|_{L^\infty}\lesssim \sum_{j'\leq j-1}2^{j'\frac{d}{q}}\big\|\dot{\Delta}_{j'} u(s)\big\|_{L^q}$~, we have that
	\begin{align}
		\left\|\dot{\Delta}_\ell~\mathcal B_t\left(\dot{S}_{j} u, \dot{\Delta}_j v\right)\right\|_{L^q}
		&\lesssim 2^\ell\sum_{j'\leq j-1}2^{j'\frac{d}{q}}\int_0^te^{-c\nu 2^{2\ell}(t-s)}\big\|\dot{\Delta}_{j'} u(s)\big\|_{L^q}\big\|\dot{\Delta}_{j} v(s)\big\|_{L^q}{\rm d}s~.\nonumber
	\end{align}
	Take $L^r$-norm in time, by Young's inequality with exponent $1+\frac{1}{r}=\frac{1}{r}+\left(\frac{1}{p}+\frac{1}{p^*}\right)$, we write
	\begin{align}
		&\left\|\dot{\Delta}_\ell~\mathcal B\left(\dot{S}_{j} u, \dot{\Delta}_j v\right)\right\|_{L^r(0,T;L^q)}\nonumber\\
		\lesssim & ~2^\ell\left(\int_0^Te^{-rc\nu 2^{2\ell}t}{\rm d}t\right)^{\frac{1}{r}}\sum_{j'\leq j-1}2^{j'\frac{d}{q}}\big\|\dot{\Delta}_{j'} u\big\|_{L^p(0,T;L^q)}\big\|\dot{\Delta}_j v\big\|_{L^{p^*}(0,T;L^q)}\nonumber\\
		\lesssim & ~2^\ell\cdot\nu^{-\frac{1}{r}}2^{-\ell\frac{2}{r}}\sum_{j'\leq j-1}2^{j'\left(\frac{d}{q}-\alpha\right)}~2^{-j\alpha'}\|u\|_{\widetilde{L}^p(0,T;\dot{B}_{q,\infty}^\alpha)}\|v\|_{\widetilde{L}^{p^*}(0,T;\dot{B}_{q,\infty}^{\alpha'})}\nonumber\\
		\lesssim & ~\nu^{-\frac{1}{r}}~2^{\ell\left(1-\frac{2}{r}\right)}2^{j\left(\frac{d}{q}-\alpha-\alpha'\right)}\|u\|_{\widetilde{L}^p(0,T;\dot{B}_{q,\infty}^\alpha)}\|v\|_{\widetilde{L}^{p^*}(0,T;\dot{B}_{q,\infty}^{\alpha'})}\nonumber
	\end{align}
	where in the last inequality we use that $\displaystyle\sum_{j'\leq j-1}2^{j'\left(\frac{d}{q}-\alpha\right)}\lesssim 2^{j\left(\frac{d}{q}-\alpha\right)}$ provided that $\alpha<\frac{d}{q}$~. Sum up in $j$, we see that
	\begin{align}
		&\sum_{j\geq\ell-N_0}\left\|\dot{\Delta}_\ell~\mathcal B\left(\dot{S}_{j} u, \dot{\Delta}_j v\right)\right\|_{L^r(0,T;L^q)}\nonumber\\
		\lesssim &~\nu^{-\frac{1}{r}}~2^{\ell\left(1-\frac{2}{r}\right)}\sum_{j\geq\ell-N}2^{j\left(\frac{d}{q}-\alpha-\alpha'\right)}\|u\|_{\widetilde{L}^p(0,T;\dot{B}_{q,\infty}^\alpha)}\|v\|_{\widetilde{L}^{p^*}(0,T;\dot{B}_{q,\infty}^{\alpha'})}\nonumber\\
		\lesssim &~\nu^{-\frac{1}{r}}~2^{\ell\left(1-\frac{2}{r}+\frac{d}{q}-\alpha-\alpha'\right)}\|u\|_{\widetilde{L}^p(0,T;\dot{B}_{q,\infty}^\alpha)}\|v\|_{\widetilde{L}^{p^*}(0,T;\dot{B}_{q,\infty}^{\alpha'})}\label{para 1}
	\end{align}
	provided $\alpha+\alpha'>\frac{d}{q}$. Estimate of $L^{r^*}(0,T; L^q)$-norm can be derive in the same way, but be awared that we need to avoid the requirement $\alpha'<\frac{d}{q}$ as we already require that $\alpha<\frac{d}{q}<\alpha+\alpha'$. For this we use Young's inequality with exponent $1+\frac{1}{r^*}=\frac{1}{r^*}+\left(\frac{1}{p}+\frac{1}{p^*}\right)$ where in the bracket the exponent $\frac{1}{p}$ is for $u$ and $\frac{1}{p^*}$ for $v$. In this way, we can derive
	\begin{align}
		\left\|\dot{\Delta}_\ell~\mathcal B\left(\dot{S}_{j} u, \dot{\Delta}_j v\right)\right\|_{L^{r^*}(0,T;L^q)}
		\lesssim \nu^{-\frac{1}{r^*}}~2^{\ell\left(1-\frac{2}{r^*}\right)}2^{j\left(\frac{d}{q}-\alpha-\alpha'\right)}\|u\|_{\widetilde{L}^p(0,T;\dot{B}_{q,\infty}^{\alpha})}\|v\|_{\widetilde{L}^{p^*}(0,T;\dot{B}_{q,\infty}^{\alpha'})}\nonumber
	\end{align}
	and then
	\begin{align}
		\sum_{j\geq\ell-N_0}\left\|\dot{\Delta}_\ell~\mathcal B\left(\dot{S}_{j} u, \dot{\Delta}_j v\right)\right\|_{L^{r^*}(0,T;L^q)}
		\lesssim\nu^{-\frac{1}{r^*}}~2^{\ell\left(1-\frac{2}{r^*}+\frac{d}{q}-\alpha-\alpha'\right)}\|u\|_{\widetilde{L}^p(0,T;\dot{B}_{q,\infty}^\alpha)}\|v\|_{\widetilde{L}^{p^*}(0,T;\dot{B}_{q,\infty}^{\alpha'})}\label{para 2}
	\end{align}
	just under the same condition $\alpha<\frac{d}{q}<\alpha+\alpha'$. To this end, one can just set
	\begin{align}
		\alpha=-1+\frac{2}{p}+\frac{d}{q}\quad and\quad \alpha'=-1+\frac{2}{p^*}+\frac{d}{q}\quad (p>2)\nonumber
	\end{align}
	to fulfil the condition. Then combining (\ref{para 1}) and (\ref{para 2}) yields that
	\begin{align}
		\left\|\mathcal Q^1(u,v)\right\|_{\Xi_{r,q}(T)}
		:&=\nu^{\frac{1}{r}}\left\|\mathcal Q^1(u,v)\right\|_{\widetilde{L}^r(0,T;\dot{B}_{q,\infty}^\alpha)}+\nu^{\frac{1}{r*}}\left\|\mathcal Q^1(u,v)\right\|_{\widetilde{L}^{r^*}(0,T;\dot{B}_{q,\infty}^{\alpha'})}\nonumber\\
		&\lesssim \|u\|_{\widetilde{L}^p(0,T;\dot{B}_{q,\infty}^\alpha)}\|v\|_{\widetilde{L}^{p^*}(0,T;\dot{B}_{q,\infty}^{\alpha'})}\nonumber\\
		&\lesssim \nu^{-1}\|u\|_{\Xi_{p,q}(T)}\|v\|_{\Xi_{p,q}(T)}~.\nonumber
	\end{align}
    
\noindent$\bf{\bullet~Weak~continuity~in~time:~}$ This is proved by the Dominated Convergence Theorem. Fix $u,v\in\Xi_{p,q}(T)$ with $p\in(2,\infty]$ and take $\mathcal Q^1(u,v)$ for example again. Then for any $\varphi\in\dot{B}_{q^*,1}^{1-d/q}$, we write by definition:
\begin{align}
    \big\langle \mathcal Q^1_t(u,v),\varphi\big\rangle=\sum_{\ell}\sum_{j\geq \ell-N_0}\left\langle \dot{\Delta}_\ell~\mathcal B_t\left(\dot{S}_{j} u, \dot{\Delta}_j v\right),\varphi\right\rangle~.\nonumber
\end{align}
We first show for each fixed $\ell,j$ that the function $t\mapsto \left\langle \dot{\Delta}_\ell~\mathcal B_t\left(\dot{S}_{j} u, \dot{\Delta}_j v\right),\varphi\right\rangle$ is continuous on $[0,T]$. We write for $0\leq t_1\leq t_2\leq T$ that
\begin{align}
    &\dot{\Delta}_\ell~\mathcal B_{t_2}\left(\dot{S}_{j} u, \dot{\Delta}_j v\right)-\dot{\Delta}_\ell~\mathcal B_{t_1}\left(\dot{S}_{j} u, \dot{\Delta}_j v\right)\nonumber\\
    =&\underbrace{\int_{t_1}^{t_2}\dot{\Delta}_\ell e^{\nu(t_2-s)\Delta}\mathcal{P}{\rm div}\left(\dot{S}_{j}u(s)\otimes\dot{\Delta}_j v(s)\right){\rm d}s}_{=:\mathcal I_{\ell,j}(t_2,t_1)}\nonumber\\
    &+\underbrace{\int_0^{t_1}\dot{\Delta}_\ell \big(e^{\nu(t_2-t_1)\Delta}-{\rm Id}\big)e^{\nu(t_1-s)\Delta}\mathcal{P}{\rm div}\left(\dot{S}_{j}u(s)\otimes\dot{\Delta}_j v(s)\right){\rm d}s}_{=:\mathcal J_{\ell,j}(t_2,t_1)}~.\nonumber
\end{align}
For the first term, by H\"older's inequality, lemma \ref{Lem. heat estimates} and \ref{Lem. homogeneous multipliers} and standard Bernstein inequality, we have:
\begin{align}
    \left|\big\langle \mathcal I_{\ell,j}(t_2,t_1),\varphi\big\rangle\right|
    &\leq\sum_{\ell'=\ell-1}^{\ell+1}\int_{t_1}^{t_2}\left|\Big\langle\dot{\Delta}_\ell e^{\nu(t_2-s)\Delta}\mathcal{P}{\rm div}\left(\dot{S}_{j}u(s)\otimes\dot{\Delta}_j v(s)\right),\dot{\Delta}_{\ell'}\varphi\Big\rangle\right|{\rm d}s\nonumber\\
    &\lesssim\sum_{\ell'=\ell-1}^{\ell+1}2^{\ell}\big\|\dot{\Delta}_{\ell'}\varphi\big\|_{L^{q^*}}\int_{t_1}^{t_2}\sum_{j'\leq j-1}2^{j'\frac{d}{q}}\big\|\dot{\Delta}_{j'} u(s)\big\|_{L^q}\big\|\dot{\Delta}_{j} v(s)\big\|_{L^q}{\rm d}s\nonumber\\
    &\lesssim\sum_{\ell'=\ell-1}^{\ell+1}2^{\ell}\big\|\dot{\Delta}_{\ell'}\varphi\big\|_{L^{q^*}}\sum_{j'\leq j-1}2^{j'\frac{d}{q}}\big\|\dot{\Delta}_{j'} u\big\|_{L^p(t_1,t_2;L^q)}\big\|\dot{\Delta}_j v\big\|_{L^{p^*}(t_1,t_2;L^q)}\nonumber\\
    &\lesssim\sum_{\ell'=\ell-1}^{\ell+1}2^{\ell}\big\|\dot{\Delta}_{\ell'}\varphi\big\|_{L^{q^*}}2^{j(1-\frac{2}{p})}\|u\|_{\Xi_{p,q}(T)}\big\|\dot{\Delta}_j v\big\|_{L^{p^*}(t_1,t_2;L^q)}\longrightarrow0\nonumber
\end{align}
when $|t_2-t_1|\rightarrow0$. Here for the last inequality we have used the fact that $p>2$ and so $\displaystyle\sum_{j'\leq j-1}2^{j'(1-\frac{2}{p})}\lesssim2^{j(1-\frac{2}{p})}$. For the other term, we have similarly that
\begin{align}
    \left|\big\langle \mathcal J_{\ell,j}(t_2,t_1),\varphi\big\rangle\right|
    &\leq\sum_{\ell'=\ell-1}^{\ell+1}\int_0^{t_1}\left|\Big\langle\dot{\Delta}_\ell e^{\nu(t_1-s)\Delta}\mathcal{P}{\rm div}\left(\dot{S}_{j}u(s)\otimes\dot{\Delta}_j v(s)\right),\big(e^{\nu(t_2-t_1)\Delta}-{\rm Id}\big)\dot{\Delta}_{\ell'}\varphi\Big\rangle\right|{\rm d}s\nonumber\\
    &\lesssim\sum_{\ell'=\ell-1}^{\ell+1}2^{\ell}\big\|\big(e^{\nu(t_2-t_1)\Delta}-{\rm Id}\big)\dot{\Delta}_{\ell'}\varphi\big\|_{L^{q^*}}2^{j(1-\frac{2}{p})}\|u\|_{\Xi_{p,q}(T)}\big\|\dot{\Delta}_j v\big\|_{L^{p^*}(0,T;L^q)}\longrightarrow0\nonumber
\end{align}
when $|t_2-t_1|\rightarrow0$. Next, we show that each $\left\langle \dot{\Delta}_\ell~\mathcal B_t\left(\dot{S}_{j} u, \dot{\Delta}_j v\right),\varphi\right\rangle$ is bounded by a summable quantity independent of $t\in[0,T]$. Again, by H\"older's inequality, lemma \ref{Lem. heat estimates} and \ref{Lem. homogeneous multipliers} and standard Bernstein inequality, we write that
\begin{align}
    \left|\left\langle \dot{\Delta}_\ell~\mathcal B_t\left(\dot{S}_{j} u, \dot{\Delta}_j v\right),\varphi\right\rangle\right|
    &\leq \sum_{\ell'=\ell-1}^{\ell+1}\int_0^t\left|\Big\langle\dot{\Delta}_\ell e^{\nu(t-s)\Delta}\mathcal{P}{\rm div}\left(\dot{S}_{j}u(s)\otimes\dot{\Delta}_j v(s)\right),\dot{\Delta}_{\ell'}\varphi\Big\rangle\right|{\rm d}s\nonumber\\
    &\lesssim\sum_{\ell'=\ell-1}^{\ell+1}2^{\ell}\big\|\dot{\Delta}_{\ell'}\varphi\big\|_{L^{q^*}}\sum_{j'\leq j-1}2^{j'\frac{d}{q}}\big\|\dot{\Delta}_{j'} u\big\|_{L^p(0,T;L^q)}\big\|\dot{\Delta}_j v\big\|_{L^{p^*}(0,T;L^q)}\nonumber\\
    &\lesssim\sum_{\ell'=\ell-1}^{\ell+1}2^{\ell}\big\|\dot{\Delta}_{\ell'}\varphi\big\|_{L^{q^*}}\sum_{j'\leq j-1}2^{j'(1-\frac{2}{p})}2^{j(1-\frac{2}{p}-\frac{d}{q})}\|u\|_{\Xi_{p,q}(T)}\|v\|_{\Xi_{p,q}(T)}\nonumber\\
    &\lesssim\sum_{\ell'=\ell-1}^{\ell+1}2^{\ell}~2^{-j\frac{d}{q}}\big\|\dot{\Delta}_{\ell'}\varphi\big\|_{L^{q^*}}\|u\|_{\Xi_{p,q}(T)}\|v\|_{\Xi_{p,q}(T)}~.\nonumber
\end{align}
The far right-hand-side is clearly summable as
\begin{align}
    \sum_{\ell\in\mathbb Z}\sum_{j\geq \ell-N_0}\sum_{\ell'=\ell-1}^{\ell+1}2^{\ell}~2^{-j\frac{d}{q}}\big\|\dot{\Delta}_{\ell'}\varphi\big\|_{L^{q^*}}\lesssim \|\varphi\|_{\dot{B}_{q^*,1}^{1-d/q}}<\infty~.\nonumber
\end{align}
Hence, Dominated Convergence Theorem implies that
\begin{align}
    \lim_{t\rightarrow t_0}\big\langle \mathcal Q^1_t(u,v),\varphi\big\rangle=\sum_{\ell}\sum_{j\geq \ell-N_0} \lim_{t\rightarrow t_0}\left\langle \dot{\Delta}_\ell~\mathcal B_t\left(\dot{S}_{j} u, \dot{\Delta}_j v\right),\varphi\right\rangle=\big\langle \mathcal Q^1_{t_0}(u,v),\varphi\big\rangle\nonumber
\end{align}
for any $t_0\in[0,T]$. This ends the proof.
\end{proof}

\
\section{Restoring well-posedness for datum in the ``bad" space $\dot{B}_{\infty,\infty}^{-1}$}
In this section, we prove Theorem \ref{Thm. restoring well-posedness by noise} and \ref{Thm. restoring well-posedness by noise II}. For convenience, let's remind the the general form (\ref{general noise form}) of the noise term: 
\begin{equation}\label{construction of noise: multiplicative}
    \left\{\begin{aligned}
    \mathcal N(t)
    &=\varepsilon_\nu\int_0^t\big(u(s)-g(s)\big){\rm dB}_s+\mathcal N^{add}(t)\\
    g(t)
    &=e^{\nu t\Delta}u_0+\int_0^te^{\nu(t-s)\Delta}{\rm d}\mathcal N^{add}(s)
    \end{aligned}\right.
\end{equation}
where $\varepsilon_\nu:=\varepsilon\nu^{1/2}$ ( $\varepsilon\in\mathbb R\setminus\{0\}$ ) and ${\rm B}$ is a one dimensional standard Brownian motions. The additive part $\mathcal N^{add}$ will be considered for the following two cases:
\begin{equation}\label{construction of mathcal N^add II}
    \left\{\begin{aligned}
        {\bf (I)}~&~\mathcal N^{add}(t)=\varepsilon_\nu\int_0^tg(s){\rm dB}_s~.~{\rm This~case~corresponds~to~linear~multiplicative~noise,}\\
        &~{\rm and~we~have}~g(t)=\mathbf M_t(\varepsilon_\nu,{\rm B})~e^{\nu t\Delta}u_0~.\\
        ~\\
        {\bf (II)}&~\mathcal N^{add}(t)=\mathcal N^{add}_{[\widetilde{\rm B},{\rm B}^j]}(t)~{\rm is~given~by}~(\ref{construction of noise II})~{\rm and}~\big\{{\rm B},\widetilde{\rm B},{\rm B}^j\big\}~{\rm is~an~independent~family}\\
        &~{\rm of~standard~Brownian~motions.}
    \end{aligned}\right.
\end{equation}
Note that the original equation (\ref{sto. NS}) with noise term given by (\ref{construction of noise: multiplicative}) can be written in the abstract form
\begin{align}\label{abstract equation II}
    u=g+\mathcal{B}^{Sto}(u,u)
\end{align}
where 
$\mathcal{B}^{Sto}$ is given by
\begin{align}
    \mathcal{B}^{Sto}_t(v_1,v_2):=-\int_0^t\exp{\left[\varepsilon_\nu({\rm B}_t-{\rm B}_s)-\frac{\varepsilon_\nu^2}{2}(t-s)\right]}e^{\nu(t-s)\Delta}\mathcal{P}{\rm div}\big(v_1(s)\otimes v_2(s)\big){\rm d}s~.\nonumber
\end{align}
This can be viewed as stochastic restoration of the original bilinear map $\mathcal B$. In order to control the exponential increment
\begin{align}
    \mathbf D_{s,t}(\varepsilon_\nu,{\rm B}):=\exp{\left[\varepsilon_\nu({\rm B}_t-{\rm B}_s)-\frac{\varepsilon_\nu^2}{2}(t-s)\right]}~,\nonumber
\end{align}
we introduce the stopping time for $M>1$:
\begin{align}\label{stopping time to control D_s,t}
    \mathcal T_1=\mathcal T_1^{\varepsilon,M}({\rm B}):=\inf\Bigg\{T>0~: \sup_{0\leq s\leq t\leq T}\exp{\left[\varepsilon_\nu({\rm B}_t-{\rm B}_s)-\frac{\varepsilon_\nu^2}{4}(t-s)\right]}>M\Bigg\}~.
\end{align}
Then the following proposition shows that, up to the stopping time, the stochastic bilinear map $\mathcal{B}^{Sto}$ is continuous on the modified solution space
\begin{align}
    X^{p,\delta}(T):=~
    &\cap\widetilde{L}^p\big(0,T;\dot{B}_{\infty,\infty}^{-1+2/p}\big)\cap\widetilde{L}^p\big(0,T;\dot{B}_{\infty,\infty}^{-1+2/p+\delta}\big)\nonumber\\
    &\cap\widetilde{L}^{p^*}\big(0,T;\dot{B}_{\infty,\infty}^{-1+2/p^*}\big)\cap\widetilde{L}^{p^*}\big(0,T;\dot{B}_{\infty,\infty}^{-1+2/p^*+\delta}\big)\nonumber
\end{align}
with norm given by
\begin{align}
    \|u\|_{X^{p,\delta}(T)}:=~
    &\nu^{1/p}\|u\|_{\widetilde{L}^p(0,T;\dot{B}_{\infty,\infty}^{-1+2/p})}
    +\nu^{1/p}\|u\|_{\widetilde{L}^{p}(0,T;\dot{B}_{\infty,\infty}^{-1+2/p+\delta})}\nonumber\\
    &+\nu^{1/p^*}\|u\|_{\widetilde{L}^{p^*}(0,T;\dot{B}_{\infty,\infty}^{-1+2/p^*})}+\nu^{1/p^*}\|u\|_{\widetilde{L}^{p^*}(0,T;\dot{B}_{\infty,\infty}^{-1+2/p^*+\delta})}~.\label{norm of space X}
\end{align}

\begin{proposition}[Stochastic Restoration of Continuity]\label{Prop. stochastic restoration of continuity}
    Let $2<p<\infty$, $0<\delta<\min\{2/p,2/p^*\}$ and $0<T\leq\infty$. Let also $\mathcal T_1=\mathcal T_1^{\varepsilon,M}({\rm B})$ be given by (\ref{stopping time to control D_s,t}) with $0<\varepsilon^2<1$. Then the bilinear map $\mathcal{B}^{Sto}:X^{p,\delta}(\mathcal T_1\wedge T)\times X^{p,\delta}(\mathcal T_1\wedge T)\rightarrow X^{p,\delta}(\mathcal T_1\wedge T)$ is continuous and
\begin{align}
    \big\|\mathcal{B}^{Sto}(v_1,v_2)\big\|_{X^{p,\delta}(\mathcal T_1\wedge T)}\leq \frac{C_2 M}{\nu\varepsilon^2}\|v_1\|_{X^{p,\delta}(\mathcal T_1\wedge T)}\|v_2\|_{X^{p,\delta}(\mathcal T_1\wedge T)}~.\nonumber
\end{align}
where $C_2>0$ is a constant possibly depends only on $p,\delta$ and $d$.
\end{proposition}

The key of the proof is that the modified norms $\|\cdot\|_{\widetilde{L}^{p^*}(0,T;\dot{B}_{\infty,\infty}^{-1+2/p^*+\delta})}$ and $\|\cdot\|_{\widetilde{L}^p(0,T;\dot{B}_{\infty,\infty}^{-1+2/p+\delta})}$ help us to control high frequency while the noise term restore the estimates on low frequency. We delay the proof of  Proposition \ref{Prop. stochastic restoration of continuity} until subsection \ref{subsec. Stochastic restoration of continuity}. With this proposition in hand, by combining with Proposition \ref{Prop. stochastic restoration of continuity} and the contraction lemma \ref{Lem. contraction}, it is not hard to see that if one can find a proper stopping time $\mathcal T$ such that
\begin{align}\label{smallness condition of g: multiplicative}
    \frac{4C_2 M}{\nu\varepsilon^2}\|g\|_{X^{p,\delta}(\mathcal T)}<1~,
\end{align}
then (\ref{abstract equation II}) has a unique solution $u\in \overline{{\rm B}\big(0,2\|g\|;~X^{p,\delta}(\mathcal T)\big)}$. Hence, to control $g$, we introduce the stopping time notation for any Brownian motion ${\rm B}$ and parameters $\kappa\neq0$, $0<\theta<1$, $A>1$:
\begin{align}\label{stopping time to control M'_t}
    \tau^{\theta,\kappa}_{A}({\rm B})=\inf\Bigg\{t>0~:~\exp\left(\kappa{\rm B}_t-\theta\frac{\kappa^2}{2}t\right)>A\Bigg\}~.
\end{align}
For case $\bf (I)$ in (\ref{construction of mathcal N^add II}), i.e. when (\ref{sto. NS}) is driven by a single linear multiplicative noise, we set
\begin{equation}\label{stopping time for g: multiplicative I}
    \tau=\tau^{1/2,\varepsilon_\nu}_{M}({\rm B})~;
\end{equation}
for case $\bf (II)$ in (\ref{construction of mathcal N^add II}), we set
\begin{equation}\label{stopping time for g: multiplicative II}
    \left\{\begin{aligned}
        \widetilde\tau
    &=\tau^{\theta,\kappa}_{A}(\widetilde{\rm B})\\
    \tau_j
    &=\tau^{1/2,\eta_j}_{A_j}({\rm B}^j)~,~j>N\\
    \tau
    &=\widetilde\tau\wedge\bigwedge_{j>N}\tau^{1/2,\eta_j}_{A_j}({\rm B}^j)~.
    \end{aligned}\right.
\end{equation}
The final stopping time is given by $\mathcal T=\mathcal T_1\wedge\tau$. To this end, one can combine the controlling technique used in the proofs of Proposition \ref{Prop. smallness of g I} and \ref{Prop. smallness of g II} to achieve the smallness condition (\ref{smallness condition of g: multiplicative}) and then estimate the probability $\mathbb P\big(\tau=\infty\big)$. To derive a lower bound for $\mathbb P\big(\mathcal T_1>T\big)$ ($0<T<\infty$), we need the following lemma.

\begin{lemma}\label{lem. probability II}
    Let $0<T<\infty$ and also $\mathcal T_1=\mathcal T_1^{\varepsilon,M}({\rm B})$ be given by (\ref{stopping time to control D_s,t}) with $\varepsilon\neq0$ and $M=e^{\varepsilon_\nu^2T/4}$. There exist universal constants $r_0,c_0>0$ such that:
    \begin{align}
        \mathbb P\big(\mathcal T_1>T\big)\geq 1-e^{-c_0\varepsilon_\nu^2T}~,\quad \varepsilon_\nu^2T\geq r_0~.\nonumber
    \end{align}
\end{lemma}

\begin{remark}
    We note that by the fact $\big\{\mathcal T_1=\infty\big\}\subset\Big\{\varepsilon_\nu({\rm B}_{n}-{\rm B}_{n-1})\leq \frac{\varepsilon_\nu^2}{4}+\ln M~,~\forall~n\in\mathbb N\Big\}$ and independence of Brownian increments, it is not hard to see that $\mathbb P\big(\mathcal T_1=\infty\big)=0$ for any fixed $M>1$ and $\kappa\neq0$.
\end{remark}

\noindent Finally, Theorem \ref{Thm. restoring well-posedness by noise} and \ref{Thm. restoring well-posedness by noise II} follow respectively by the two propositions below. For consistency of notations, we would denote by $C=C(p,d)>0$ some constant and ${\bf R}_\infty:=\nu^{-1}\|u_0\|_{\dot{B}_{\infty,\infty}^{-1}}$.

\begin{proposition}\label{Prop. smallness of g: multiplicative I}
    Let $2<p<\infty$ and $u_0\in \dot{B}_{\infty,\infty}^{-1+\delta}$ with $0<\delta<\min\{2/p,2/p^*\}$. Let the noise term be given by (\ref{construction of noise: multiplicative}) and case (I) of (\ref{construction of mathcal N^add II}) with $0<\varepsilon^2<1$. Also, let $\mathcal T_1$ be given by (\ref{stopping time to control D_s,t}) and  $\tau$ given by (\ref{stopping time for g: multiplicative I}) with $M>1$. Then up the stopping time $\mathcal T=\mathcal T_1\wedge\tau$, we have
    \begin{align}\label{smallness of g: multiplicative I}
        \|g\|_{X^{p,\delta}(\mathcal T)}\lesssim_{p,\delta,d} \frac{M}{\varepsilon^2}\|u_0\|_{\dot{B}_{\infty,\infty}^{-1+\delta}}~.
    \end{align}
    Moreover, if we set $M=e^{\varepsilon_\nu^2T/4}$, then there exists universal constants $r_0',c_0'>0$ such that
    \begin{align}
        \mathbb P\big(\mathcal T>T\big)\geq 1-2e^{-c_0'\varepsilon_\nu^2T}~.\nonumber
    \end{align}
    holds for all finite $T\geq r_0'\varepsilon_\nu^{-2}$.
\end{proposition}

\begin{proposition}\label{Prop. smallness of g: multiplicative II}
    Let $2<p<\infty$, $0<\delta<\min\{2/p,2/p^*\}$ and $u_0\in \dot{B}_{\infty,\infty}^{-1}$. Let the noise term be given by (\ref{construction of noise: multiplicative}) and case (II) of (\ref{construction of mathcal N^add II}); the stopping time $\mathcal T=\mathcal T_1\wedge\tau$ given by (\ref{stopping time to control D_s,t}) and (\ref{stopping time for g: multiplicative II}), with parameters $N\in\mathbb N$, $\varepsilon^2,\theta\in(0,1)$, $M>1$,
    \begin{align}
        \kappa^2
        &>\frac{C\nu 2^{4N+1}}{1-\theta}\left(\frac{M}{\varepsilon^2}\right)^{p+p^*}\left({\bf R}_\infty^p+{\bf R}_\infty^{p^*}\right)=:K_\infty^2~,\nonumber\\
        \eta_j
        &=C\left(\frac{M}{\varepsilon^2}\right)^{\frac{p+p^*}{2}}\sqrt{{\bf R}_\infty^p+{\bf R}_\infty^{p^*}}~\nu^{1/2} 2^{j(1+\gamma)}~,\quad \gamma>\frac{\delta}{2}(p+p^*)~,\nonumber
    \end{align}
    and $A=(\kappa/K_\infty)^{\frac{2}{p+p^*}}$, $A_j=\left(\frac{\pi^2}{12}j^2\right)^2$~. Then there exist integer constant $N_0=N_0(p,\gamma,\delta)>0$ such that for any $N\geq N_0$, we have
    \begin{align}\label{smallness condition of g: B^{-1} II}
        \frac{4C_2 M}{\nu\varepsilon^2}\|g\|_{X^{p,\delta}(\mathcal T)}<1~.
    \end{align}
    Moreover, if we set $M=e^{\varepsilon_\nu^2T/4}$, then there exists universal constants $r_0,c_0>0$ such that
    \begin{align}\label{probability lower bound: B^{-1} II}
        \mathbb P\big(\mathcal T>T\big)\geq \left(1-e^{-c_0\varepsilon_\nu^2T}\right)\left[1-\left(\frac{\kappa}{K_\infty}\right)^{-\frac{2\theta}{p+p^*}}\right]\exp\Bigg(-\frac{6}{\pi^2}\sum_{j>N}j^{-2}\Bigg)
    \end{align}
    holds for all finite $T\geq r_0\varepsilon_\nu^{-2}$.
\end{proposition}

\noindent Finally, the following proposition indicates that our solution has modified paths $u-g\in C_{w^*,loc}\big([0,\mathcal T);\dot{B}_{\infty,\infty}^{-1+\delta}\big)$.

\begin{proposition}\label{Prop. stochastic restoration of continuity II}
    Let the exponents $p,\delta,\varepsilon$ and also the stopping time $\mathcal T_1$ be given as in Proposition \ref{Prop. stochastic restoration of continuity}. Set $\mathcal T_1:=\mathcal T_1\wedge T$ with arbitrarily fixed $0<T\leq\infty$. The bilinear map $\mathcal{B}^{Sto}:X^{p,\delta}(\mathcal T_1')\times X^{p,\delta}(\mathcal T_1')\rightarrow {L}^{\infty}\big(0,\mathcal T_1';\dot{B}_{\infty,\infty}^{-1+\delta}\big)\cap {L}^{1}\big(0,\mathcal T_1';\dot{B}_{\infty,\infty}^{1}\big)$ is continuous and
\begin{align}
    \big\|\mathcal{B}^{Sto}(v_1,v_2)\big\|_{{L}^{\infty}(0,\mathcal T_1';\dot{B}_{\infty,\infty}^{-1+\delta})}+\varepsilon^2\nu\big\|\mathcal{B}^{Sto}(v_1,v_2)\big\|_{{L}^{1}(0,\mathcal T_1';\dot{B}_{\infty,\infty}^{1})}\leq \frac{C_2 M}{\nu}\|v_1\|_{X^{p,\delta}(\mathcal T_1')}\|v_2\|_{X^{p,\delta}(\mathcal T_1')}~.\nonumber
\end{align}
Furthermore, for each $\varphi\in\dot{B}_{1,1}^{1-\delta}$, the function $t\mapsto \big\langle B^{Sto}_t(v_1,v_2),\varphi\big\rangle$ is continuous on $[0,\mathcal T_1']$ provided that $v_1,v_2\in X^{p,\delta}(\mathcal T_1')$.
\end{proposition}\

\par For the organization of the rest of the section, we first prove Proposition \ref{Prop. stochastic restoration of continuity} and  \ref{Prop. stochastic restoration of continuity II} in the up-coming subsection \ref{subsec. Stochastic restoration of continuity}; and then in subsection \ref{subsec. Domination by noise II} we prove Proposition \ref{Prop. smallness of g: multiplicative I} and \ref{Prop. smallness of g: multiplicative II}; finally we present the proof of lemma \ref{lem. probability II} in subsection \ref{subsec. The probability lemma}.

\
\subsection{Stochastic restoration of continuity}\label{subsec. Stochastic restoration of continuity}
\begin{proof}[Proof of Proposition \ref{Prop. stochastic restoration of continuity}]
    Let $2<p<\infty$, $0<\delta<\min\{2/p,2/p^*\}$, $0<T\leq\infty$ and $0<\varepsilon^2<1$. Let also the stopping time $\mathcal T_1=\mathcal T_1^{\kappa,M}({\rm B})$ be given by (\ref{stopping time to control D_s,t}) with $M>1$. To keep notation simple, we denote $\mathcal T_1'=\mathcal T_1\wedge T$. We recall the Bony decomposition (\ref{Bony decomposition}) which we now apply to the bilinear map $\mathcal{B}^{Sto}$. As in the proof of Proposition \ref{Prop. continuity of bilinear map}, we only show the estimates of the term 
    \begin{align}
        \mathcal{Q}^{Sto.1}(v_1,v_2):=\sum_{\ell\in\mathbb Z}\sum_{j\geq\ell-N_0}\dot{\Delta}_\ell~\mathcal{B}^{Sto}\left(\dot{S}_{j} v_1, \dot{\Delta}_j v_1\right)~.\nonumber
    \end{align}
    For $0\leq s\leq t\leq\mathcal T_1'$, we obviously have
    $\mathbf D_{s,t}(\varepsilon_\nu,{\rm B})\leq Me^{-\frac{\varepsilon^2\nu}{4}(t-s)}$. Then together by lemma \ref{Lem. heat estimates}, we have that
    \begin{align}
        \left\|\dot{\Delta}_\ell~\mathcal{B}^{Sto}_t\left(\dot{S}_{j} v_1, \dot{\Delta}_j v_1\right)\right\|_{L^\infty}\lesssim M2^\ell\int_0^te^{-(t-s)\left(\frac{\varepsilon^2\nu}{4}+c\nu2^{2\ell}\right)}\sum_{j'\leq j-1}\big\|\dot{\Delta}_{j'}v_1(s)\big\|_{L^\infty}\big\|\dot{\Delta}_{j}v_2(s)\big\|_{L^\infty}{\rm d}s~.\nonumber
    \end{align}
    Take $L^p(0,\mathcal T_1')$-norm on both side and apply Young's inequality with exponents $1+\frac{1}{p}=\frac{1}{p}+\left(\frac{1}{p}+\frac{1}{p^*}\right)$, we see that
    \begin{align}
        &\left\|\dot{\Delta}_\ell~\mathcal{B}^{Sto}\left(\dot{S}_{j} v_1, \dot{\Delta}_j v_1\right)\right\|_{L^p(0,\mathcal T_1'; L^\infty)}\nonumber\\
        \lesssim &~ M2^\ell\left(\frac{\varepsilon^2\nu}{4}+c\nu2^{2\ell}\right)^{-\frac{1}{p}}\sum_{j'\leq j-1}2^{j'\left(1-\frac{2}{p}\right)}\|v_1\|_{\widetilde{L}^p(0,\mathcal T_1';\dot{B}_{\infty,\infty}^{-1+2/p})}2^{j\left(1-\frac{2}{p^*}-\delta\right)}\|v_2\|_{\widetilde{L}^{p^*}(0,\mathcal T_1';\dot{B}_{\infty,\infty}^{-1+2/p^*+\delta})}\nonumber\\
        \lesssim &~ M2^\ell\left(\frac{\varepsilon^2\nu}{4}+c\nu2^{2\ell}\right)^{-\frac{1}{p}}2^{-j\delta}\|v_1\|_{\widetilde{L}^p(0,\mathcal T_1';\dot{B}_{\infty,\infty}^{-1+2/p})}\|v_2\|_{\widetilde{L}^{p^*}(0,\mathcal T_1';\dot{B}_{\infty,\infty}^{-1+2/p^*+\delta})}~.\nonumber
    \end{align}
    Here we have used the fact that $1-\frac{2}{p}>0$. Sum in $j$ then we get
    \begin{align}
        &2^{\ell\left(-1+\frac{2}{p}\right)}\left\|\sum_{j\geq\ell-N_0}\dot{\Delta}_\ell~\mathcal{B}^{Sto}\left(\dot{S}_{j} v_1, \dot{\Delta}_j v_1\right)\right\|_{L^p(0,\mathcal T_1'; L^\infty)}\nonumber\\
        \lesssim &~\frac{M2^{\ell\left(\frac{2}{p}-\delta\right)}}{\left(\frac{\varepsilon^2\nu}{4}+c\nu2^{2\ell}\right)^{\frac{1}{p}}}\|v_1\|_{\widetilde{L}^p(0,\mathcal T_1';\dot{B}_{\infty,\infty}^{-1+2/p})}\|v_2\|_{\widetilde{L}^{p^*}(0,\mathcal T_1';\dot{B}_{\infty,\infty}^{-1+2/p^*+\delta})}~.\nonumber
    \end{align}
    As $\frac{2}{p}-\delta>0$, the right-hand-side is obviously bounded by constant multiple of 
    \begin{align}\label{dominate quantity of B^Sto}
        \frac{M}{\big(\nu\varepsilon^2\big)^{\frac{1}{p}}}\|v_1\|_{\widetilde{L}^p(0,\mathcal T_1';\dot{B}_{\infty,\infty}^{-1+2/p})}\|v_2\|_{\widetilde{L}^{p^*}(0,\mathcal T_1';\dot{B}_{\infty,\infty}^{-1+2/p^*+\delta})}
    \end{align} 
    when $\ell\leq 0$. For $\ell>0$, we use the fact that $0<\varepsilon^2<1$ so that $\frac{\varepsilon^2\nu}{4}+c\nu2^{2\ell}>c\varepsilon^2\nu2^{2\ell}$, then one sees that the right-hand-side is again bounded by the same quantity. Hence, again use the fact that $0<\varepsilon^2<1$, we get
    \begin{align}
        \nu^{1/p}\left\|\mathcal{Q}^{Sto.1}(v_1,v_2)\right\|_{\widetilde{L}^p(0,\mathcal T_1';\dot{B}_{\infty,\infty}^{-1+2/p})}
        &\lesssim \frac{M}{\varepsilon^2}\|v_1\|_{\widetilde{L}^p(0,\mathcal T_1';\dot{B}_{\infty,\infty}^{-1+2/p})}\|v_2\|_{\widetilde{L}^{p^*}(0,\mathcal T_1';\dot{B}_{\infty,\infty}^{-1+2/p^*+\delta})}\nonumber\\
        &\lesssim \frac{M}{\nu\varepsilon^2}\|v_1\|_{X^{p,\delta}(\mathcal T_1')}\|v_2\|_{X^{p,\delta}(\mathcal T_1')}~.\nonumber
    \end{align}
    Similarly, one can apply Young's inequality with exponents $1+\frac{1}{p^*}=\frac{1}{p^*}+\left(\frac{1}{p}+\frac{1}{p^*}\right)$ instead to show that
    \begin{align}
        &2^{\ell\left(-1+\frac{2}{p^*}+\delta\right)}\left\|\sum_{j\geq\ell-N_0}\dot{\Delta}_\ell~\mathcal{B}^{Sto}\left(\dot{S}_{j} v_1, \dot{\Delta}_j v_1\right)\right\|_{L^{p^*}(0,\mathcal T_1'; L^\infty)}\nonumber\\
        \lesssim &~\frac{M2^{\ell\frac{2}{p^*}}}{\left(\frac{\varepsilon^2\nu}{4}+c\nu2^{2\ell}\right)^{\frac{1}{p^*}}}\|v_1\|_{\widetilde{L}^p(0,\mathcal T_1';\dot{B}_{\infty,\infty}^{-1+2/p})}\|v_2\|_{\widetilde{L}^{p^*}(0,\mathcal T_1';\dot{B}_{\infty,\infty}^{-1+2/p^*+\delta})}~.\nonumber
    \end{align}
    Then for the same reason, the right-hand-side is bounded by constant multiple of the quantity (\ref{dominate quantity of B^Sto}) with $p$ being replaced by $p^*$, and so
    \begin{align}
        \nu^{1/p^*}\left\|\mathcal{Q}^{Sto.1}(v_1,v_2)\right\|_{\widetilde{L}^{p^*}(0,\mathcal T_1';\dot{B}_{\infty,\infty}^{-1+2/p^*+\delta})}
        \lesssim \frac{M}{\nu\varepsilon^2}\|v_1\|_{X^{p,\delta}(\mathcal T_1')}\|v_2\|_{X^{p,\delta}(\mathcal T_1')}~.\nonumber
    \end{align}
    Estimates of the norms $\nu^{1/p}\|\cdot\|_{\widetilde{L}^{p}(0,\mathcal T_1';\dot{B}_{\infty,\infty}^{-1+2/p+\delta})}$ and $\nu^{1/p^*}\|\cdot\|_{\widetilde{L}^{p^*}(0,\mathcal T_1';\dot{B}_{\infty,\infty}^{-1+2/p^*})}$
    can be derived in the same way. This ends the proof.

\end{proof}

\

\begin{proof}[Proof of Proposition \ref{Prop. stochastic restoration of continuity II}]
    We use the same notations as in the previous proof. So if we take $L^\infty(0,\mathcal T_1')$-norm on both side of
    \begin{align}
        \left\|\dot{\Delta}_\ell~\mathcal{B}^{Sto}_t\left(\dot{S}_{j} v_1, \dot{\Delta}_j v_1\right)\right\|_{L^\infty}\lesssim M2^\ell\int_0^te^{-(t-s)\left(\frac{\varepsilon^2\nu}{4}+c\nu2^{2\ell}\right)}\sum_{j'\leq j-1}\big\|\dot{\Delta}_{j'}v_1(s)\big\|_{L^\infty}\big\|\dot{\Delta}_{j}v_2(s)\big\|_{L^\infty}{\rm d}s\nonumber
    \end{align}
    and apply Young's inequality with exponents $1+\frac{1}{\infty}=\frac{1}{\infty}+\left(\frac{1}{p}+\frac{1}{p^*}\right)$, we get
    \begin{align}
        &\left\|\dot{\Delta}_\ell~\mathcal{B}^{Sto}\left(\dot{S}_{j} v_1, \dot{\Delta}_j v_1\right)\right\|_{L^\infty(0,\mathcal T_1'; L^\infty)}\nonumber\\
        \lesssim &~ M2^\ell\sum_{j'\leq j-1}2^{j'\left(1-\frac{2}{p}\right)}\|v_1\|_{\widetilde{L}^p(0,\mathcal T_1';\dot{B}_{\infty,\infty}^{-1+2/p})}2^{j\left(1-2/p^*-\delta\right)}\|v_2\|_{\widetilde{L}^{p^*}(0,\mathcal T_1';\dot{B}_{\infty,\infty}^{-1+\frac{2}{p^*}+\delta})}\nonumber\\
        \lesssim &~ M2^\ell2^{-j\delta}\|v_1\|_{\widetilde{L}^p(0,\mathcal T_1';\dot{B}_{\infty,\infty}^{-1+2/p})}\|v_2\|_{\widetilde{L}^{p^*}(0,\mathcal T_1';\dot{B}_{\infty,\infty}^{-1+2/p^*+\delta})}~.\nonumber
    \end{align}
    Sum in $j$, we get
    \begin{align}
        2^{\ell(-1+\delta)}\left\|\sum_{j\geq\ell-N_0}\dot{\Delta}_\ell~\mathcal{B}^{Sto}\left(\dot{S}_{j} v_1, \dot{\Delta}_j v_1\right)\right\|_{L^\infty(0,\mathcal T_1'; L^\infty)}\lesssim M\|v_1\|_{\widetilde{L}^p(0,\mathcal T_1';\dot{B}_{\infty,\infty}^{-1+2/p})}\|v_2\|_{\widetilde{L}^{p^*}(0,\mathcal T_1';\dot{B}_{\infty,\infty}^{-1+2/p^*+\delta})}\nonumber
    \end{align}
    and  so
    \begin{align}
        \left\|\mathcal{Q}^{Sto.1}(v_1,v_2)\right\|_{{L}^{\infty}(0,\mathcal T_1';\dot{B}_{\infty,\infty}^{-1+\delta})}
        \lesssim M\nu^{-1}\|v_1\|_{X^{p,\delta}(\mathcal T_1')}\|v_2\|_{X^{p,\delta}(\mathcal T_1')}~.\nonumber
    \end{align}
    If we apply Young's inequality with exponents $1+\frac{1}{1}=\frac{1}{1}+\left(\frac{1}{p}+\frac{1}{p^*}\right)$ instead, we would get
    \begin{align}
        &\left\|\dot{\Delta}_\ell~\mathcal{B}^{Sto}\left(\dot{S}_{j} v_1, \dot{\Delta}_j v_1\right)\right\|_{L^1(0,\mathcal T_1'; L^\infty)}\nonumber\\
        \lesssim &~ M2^\ell\left(\frac{\varepsilon^2\nu}{4}+c\nu2^{2\ell}\right)^{-1}\sum_{j'\leq j-1}2^{j'\left(1-\frac{2}{p}\right)}\|v_1\|_{\widetilde{L}^p(0,\mathcal T_1';\dot{B}_{\infty,\infty}^{-1+2/p})}2^{j\left(1-2/p^*-\delta\right)}\|v_2\|_{\widetilde{L}^{p^*}(0,\mathcal T_1';\dot{B}_{\infty,\infty}^{-1+2/p^*+\delta})}\nonumber\\
        \lesssim &~ M2^\ell\left(\frac{\varepsilon^2\nu}{4}+c\nu2^{2\ell}\right)^{-1}2^{-j\delta}\|v_1\|_{\widetilde{L}^p(0,\mathcal T_1';\dot{B}_{\infty,\infty}^{-1+2/p})}\|v_2\|_{\widetilde{L}^{p^*}(0,\mathcal T_1';\dot{B}_{\infty,\infty}^{-1+2/p^*+\delta})}\nonumber
    \end{align}
    and so summing in $j$ would give 
    \begin{align}
        &2^{\ell}\left\|\sum_{j\geq\ell-N_0}\dot{\Delta}_\ell~\mathcal{B}^{Sto}\left(\dot{S}_{j} v_1, \dot{\Delta}_j v_1\right)\right\|_{L^1(0,\mathcal T_1'; L^\infty)}\nonumber\\
        \lesssim &~\frac{M2^{\ell\left(2-\delta\right)}}{\frac{\varepsilon^2\nu}{4}+c\nu2^{2\ell}}\|v_1\|_{\widetilde{L}^p(0,\mathcal T_1';\dot{B}_{\infty,\infty}^{-1+2/p})}\|v_2\|_{\widetilde{L}^{p^*}(0,\mathcal T_1';\dot{B}_{\infty,\infty}^{-1+2/p^*+\delta})}~.\nonumber
    \end{align}
    Note that $2-\delta>0$. Then by the same argument in the previous proof, we can derive
    \begin{align}
        \nu\left\|\mathcal{Q}^{Sto.1}(v_1,v_2)\right\|_{\widetilde{L}^{1}(0,\mathcal T_1';\dot{B}_{\infty,\infty}^{1})}
        \lesssim \frac{M}{\nu\varepsilon^2}\|v_1\|_{X^{p,\delta}(\mathcal T_1')}\|v_2\|_{X^{p,\delta}(\mathcal T_1')}~.\nonumber
    \end{align}
    The proof of weak continuity of $t\mapsto \mathcal{B}^{Sto}_t(v_1,v_2)$ is very similar to the weak continuity part in the proof of Proposition \ref{Prop. continuity of bilinear map}, so we just omit it here.
\end{proof}

\
\subsection{Domination by noise}\label{subsec. Domination by noise II}

\begin{proof}[Proof of Proposition \ref{Prop. smallness of g: multiplicative I}]
    Let $u_0$, the noise term and the stopping time $\mathcal T$ be given as in the proposition. We remind that $g(t)=\mathbf M_t(\varepsilon_\nu,{\rm B})~e^{\nu t\Delta}u_0$. Before the stopping time, we obviously have $\mathbf M_t(\varepsilon_\nu,{\rm B})\leq Me^{-\frac{\varepsilon^2\nu}{4}t}$. Then by lemma \ref{Lem. heat estimates},
    \begin{align}
        \left\|\mathbf M_t(\varepsilon_\nu,{\rm B})~e^{\nu t\Delta}\dot{\Delta}_ju_0\right\|_{L^\infty}\lesssim e^{-t\left(\frac{\varepsilon^2\nu}{4}+c\nu2^{2j}\right)}\|\dot{\Delta}_ju_0\|_{L^\infty}~.\nonumber
    \end{align}
    Take $L^p(0,\mathcal T)$-norm ($1\leq p<\infty$), we have that
    \begin{align}
        2^{j\left(-1+\frac{2}{p}\right)}\left\|\mathbf M(\varepsilon_\nu,{\rm B})~e^{\nu t\Delta}\dot{\Delta}_ju_0\right\|_{L^p(0,\mathcal T; L^\infty)}
        &\lesssim \frac{M2^{j\left(\frac{2}{p}-\delta\right)}}{\left(\frac{\varepsilon^2\nu}{4}+c\nu2^{2j}\right)^{\frac{1}{p}}}\|u_0\|_{\dot{B}_{\infty,\infty}^{-1+\delta}}~,\nonumber\\
        2^{j\left(-1+\frac{2}{p}+\delta\right)}\left\|\mathbf M(\varepsilon_\nu,{\rm B})~e^{\nu t\Delta}\dot{\Delta}_ju_0\right\|_{L^p(0,\mathcal T; L^\infty)}
        &\lesssim \frac{M2^{j\frac{2}{p}}}{\left(\frac{\varepsilon^2\nu}{4}+c\nu2^{2j}\right)^{\frac{1}{p}}}\|u_0\|_{\dot{B}_{\infty,\infty}^{-1+\delta}}~.\nonumber
    \end{align}
    So by the same argument in the proof of Proposition \ref{Prop. stochastic restoration of continuity}, one sees that 
    \begin{align}
        \nu^{1/p}\|g\|_{\widetilde{L}^p(0,\mathcal T;\dot{B}_{\infty,\infty}^{-1+2/p})}+\nu^{1/p}\|g\|_{\widetilde{L}^p(0,\mathcal T;\dot{B}_{\infty,\infty}^{-1+2/p})}\lesssim \frac{M}{\varepsilon^2}\|u_0\|_{\dot{B}_{\infty,\infty}^{-1+\delta}}\nonumber
    \end{align}
    and hence (\ref{smallness of g: multiplicative I}) hold true. By the way, one clearly have that $\|g\|_{L^\infty(0,\mathcal T;\dot{B}_{\infty,\infty}^{-1+\delta})}\lesssim M\|u_0\|_{\dot{B}_{\infty,\infty}^{-1+\delta}}$.
    
    Next if we set $M=e^{\varepsilon_\nu^2T/4}$, then by lemma \ref{lem. probability II}, there exist universal constants $r_0,c_0>0$ such that when $\varepsilon_\nu^2T\geq r_0$~, we have 
    \begin{align}
        \mathbb P\big(\mathcal T_1>T\big)\geq 1-e^{-c_0\varepsilon_\nu^2T}~.\nonumber
    \end{align} 
    On the other hand, we also have by lemma \ref{Lem. probability I} that
    \begin{align}
        \mathbb P\big(\tau=\infty\big)\geq 1-M^{-\frac{1}{2}}\geq 1-e^{-\frac{1}{8}\varepsilon_\nu^2T}~.\nonumber
    \end{align}
    Obviously, one can set $c_0'=\min\{c_0,1/8\}$ and choose $r_0'\geq r_0$ such that $e^{-c_0'\varepsilon_\nu^2T}\leq1/2$ whenever $\varepsilon_\nu^2T\geq r_0'$. And so
    \begin{align}
       \mathbb P\big(\mathcal T>T\big)
       &\geq \mathbb P\big(\mathcal T_1>T~and~\tau=\infty\big)\nonumber\\
       &\geq 1-\Big[1-\mathbb P\big(\mathcal T_1>T\big)\Big]-\Big[1-\mathbb P\big(\tau=\infty\big)\Big]\nonumber\\
       &\geq 1-2e^{-c_0'\varepsilon_\nu^2T}~.\nonumber
    \end{align}
\end{proof}

\begin{proof}[Proof of Proposition \ref{Prop. smallness of g: multiplicative II}]
    Let $2<p<\infty$, $0<\delta<\min\{2/p,2/p^*\}$ and $u_0\in \dot{B}_{\infty,\infty}^{-1}$. Let's briefly recall the setting for convenience. Firstly, the noise term is given by (\ref{construction of noise: multiplicative}) and case (II) of (\ref{construction of mathcal N^add II}), that is to say ($\varepsilon_\nu=\varepsilon\nu^{1/2}$) :
    \begin{align}
        \mathcal N(t)
        &=\varepsilon_\nu\int_0^t\big[u(s)-g(s)\big]{\rm dB}_s+\mathcal N^{add}(t)~,\nonumber\\
        g(t)
        &=e^{\nu t\Delta}u_0+\int_0^te^{\nu(t-s)\Delta}{\rm d}\mathcal N^{add}(s)\nonumber
    \end{align}
    with $\mathcal N^{add}(t)=\mathcal N^{add}_{[\widetilde{\rm B},{\rm B}^j]}(t)$ given by (\ref{construction of noise II}). The Brownian motions are all independent of each other. By It\^o's formula, it is not hard to verify that
    \begin{align}
        g(t)=\underbrace{\sum_{j\leq N}\mathbf M_t(\kappa,\widetilde{\rm B})~e^{\nu t\Delta}\dot\Delta_j u_0}_{={\mathcal L}^N}+\underbrace{\sum_{j> N}\mathbf M_t(\eta_j,{\rm B}^j)e^{\nu t\Delta}\dot\Delta_j u_0}_{={\mathcal H}^N}~.\nonumber
    \end{align}
    We also remind that the stopping times $\mathcal T_1=\mathcal T_1^{\varepsilon,M}({\rm B})$, $\widetilde\tau=\tau^{\theta,\kappa}_{A}(\widetilde{\rm B})$ and $\tau_j=\tau^{1/2,\eta_j}_{A_j}({\rm B}^j)$ are given in notation of (\ref{stopping time to control D_s,t}) and  (\ref{stopping time to control M'_t}) respectively, and finally
    \begin{align}
        \mathcal T=\mathcal T_1\wedge\underbrace{\widetilde\tau\wedge\bigwedge_{j>N}\tau_j}_{=\tau}~.\nonumber
    \end{align}
    The parameters $\varepsilon^2,\theta\in(0,1)$ are arbitrarily fixed; $N\in\mathbb N$, $\kappa,\eta_j\neq0$ and $A,A_j>1$ would be determined later when we control the norm $\|g\|_{X^{p,\delta}}$~; $M=e^{\varepsilon_\nu^2T/4}>1$ is determined by a further large time parameter $T>0$ which would be determined when we estimate the probability $\mathbb P\big(\mathcal T>T\big)$.\\
    
    \noindent $\mathbf{\bullet~Proof~of~(\ref{smallness condition of g: B^{-1} II}).}$ The control part is similar to the proof of Proposition \ref{Prop. smallness of g II}. Like before, we denote by $C=C(p,d)>0$ some constant that may be different from line to line. We first control the high frequency part ${\mathcal H}^N$. For this, we set $A_j=\left(\frac{\pi^2}{12}j^2\right)^2$ as before. Then up to the stopping time $\mathcal T\leq\tau_j$, we have for all $j>N$ that:
    \begin{align}
        \frac{4C_2M}{\nu\varepsilon^2}\nu^{\frac{1}{p}}2^{j\left(-1+\frac{2}{p}\right)}\big\|\mathbf M\left(\eta_j,{\rm B}^j\right)e^{\nu t\Delta}\dot\Delta_j u_0\big\|_{L^p(0,\mathcal T; L^\infty)}
        &\lesssim\left(\frac{\nu 2^{2j}}{\eta_j^2}\right)^{\frac{1}{p}}A_j~\frac{M\|u_0\|_{\dot{B}_{\infty,\infty}^{-1}}}{\nu\varepsilon^2}~,\nonumber\\
        \frac{4C_2M}{\nu\varepsilon^2}\nu^{\frac{1}{p}}2^{j\left(-1+\frac{2}{p}+\delta\right)}\big\|\mathbf M\left(\eta_j,{\rm B}^j\right)e^{\nu t\Delta}\dot\Delta_j u_0\big\|_{L^p(0,\mathcal T; L^\infty)}
        &\lesssim\left(\frac{\nu 2^{(2+p\delta)j}}{\eta_j^2}\right)^{\frac{1}{p}}A_j~\frac{M\|u_0\|_{\dot{B}_{\infty,\infty}^{-1}}}{\nu\varepsilon^2}~.\nonumber
    \end{align}
    Here we should at least set $N$ such that $A_N>1$. Note that $p\delta<2$, so if we take non-decreasing $\eta_j$ such that
	\begin{align}\label{eta_j condition: multiplicative}
		\eta_{j-1}^2\geq \left(\frac{M}{\varepsilon^2}\right)^p(9C{\bf R}_\infty)^p\nu 2^{(2+p\delta)j}\left(\frac{\pi^2}{12}(j+1)^2\right)^{2p},\quad \bigg({\bf R}_\infty=\frac{\|u_0\|_{\dot{B}_{\infty,\infty}^{-1}}}{\nu}\bigg)
	\end{align}
	the above inequalities would imply that 
	\begin{align}
		\frac{4C_2M}{\nu\varepsilon^2}\nu^{\frac{1}{p}}2^{j\left(-1+\frac{2}{p}\right)}\big\|\dot\Delta_j {\mathcal H}^N\big\|_{L^p(0,\mathcal T; L^\infty)}
		&\lesssim\sum_{\ell=j-1}^{j+1}\left(\frac{\nu 2^{2\ell}}{\eta_\ell^2}\right)^{\frac{1}{p}}A_\ell\frac{M\|u_0\|_{\dot{B}_{\infty,\infty}^{-1}}}{\nu\varepsilon^2}\nonumber\\
		&\leq C\left(\frac{\nu 2^{2j}}{\eta_{j-1}^2}\right)^{\frac{1}{p}}A_{j+1}\frac{M\|u_0\|_{\dot{B}_{\infty,\infty}^{-1}}}{\nu\varepsilon^2}\leq 1/9~,\nonumber\\
        \frac{4C_2M}{\nu\varepsilon^2}\nu^{\frac{1}{p}}2^{j\left(-1+\frac{2}{p}+\delta\right)}\big\|\dot\Delta_j {\mathcal H}^N\big\|_{L^p(0,\mathcal T; L^\infty)}
		&\lesssim\sum_{\ell=j-1}^{j+1}\left(\frac{\nu 2^{(2+p\delta)\ell}}{\eta_\ell^2}\right)^{\frac{1}{p}}A_\ell\frac{M\|u_0\|_{\dot{B}_{\infty,\infty}^{-1}}}{\nu\varepsilon^2}\nonumber\\
		&\leq C\left(\frac{\nu 2^{(2+p\delta)j}}{\eta_{j-1}^2}\right)^{\frac{1}{p}}A_{j+1}\frac{M\|u_0\|_{\dot{B}_{\infty,\infty}^{-1}}}{\nu\varepsilon^2}\leq 1/9~,\nonumber
	\end{align}
	and hence
	\begin{align}
		\frac{4C_2M}{\nu\varepsilon^2}\nu^{\frac{1}{p}}\big\|{\mathcal H}^N\big\|_{\widetilde{L}^p\big(0,\mathcal T; \dot{B}_{\infty,\infty}^{-1+2/p}\big)}+\frac{4C_2M}{\nu\varepsilon^2}\nu^{\frac{1}{p}}\big\|{\mathcal H}^N\big\|_{\widetilde{L}^p\big(0,\mathcal T; \dot{B}_{\infty,\infty}^{-1+2/p+\delta}\big)}\leq 2/9~.\nonumber
	\end{align}
    The same thing also holds if one replace $p$ by $p^*$. So in view of (\ref{eta_j condition: multiplicative}), one may just set
    \begin{align}
        \eta_j=C\left(\frac{M}{\varepsilon^2}\right)^{\frac{p+p^*}{2}}\sqrt{{\bf R}_\infty^p+{\bf R}_\infty^{p^*}}\nu^{1/2}2^{(1+\gamma)j}~,\quad \gamma>\frac{\delta}{2}(p+p^*)~.\nonumber
    \end{align}
    Then for $N\in\mathbb N$ such that
    \begin{align}
        2^{(1+\gamma)j}\geq\left(\frac{\pi^2}{12}(j+1)^2\right)^{p+p^*}2^{\frac{\delta}{2}(p+p^*)j}~,\quad \forall j>N,\nonumber
    \end{align}
    we would have $\displaystyle\frac{4C_2M}{\nu\varepsilon^2}\big\|{\mathcal H}^N\big\|_{X^{p,\delta}(\mathcal T)}\leq 4/9$~. Now fix any such $N\in\mathbb N$, we control the low frequency part ${\mathcal L}^N$. Similar as in the proof of Proposition \ref{Prop. smallness of g I}, up to the stopping time $\mathcal T\leq\widetilde\tau$, we have for all $j\leq N$ that
    \begin{align}
       \frac{4C_2M}{\nu\varepsilon^2}\nu^{\frac{1}{p}}2^{j\left(-1+\frac{2}{p}\right)}\big\|\dot\Delta_j \mathcal L^N\big\|_{L^p(0,\mathcal T; L^\infty)}
       &\leq C\left[\frac{\nu 2^{2j}}{(1-\theta)\frac{\kappa^2}{2}}\right]^{\frac{1}{p}}A~\frac{M\|u_0\|_{\dot{B}_{\infty,\infty}^{-1}}}{\nu\varepsilon^2}~,\nonumber\\
       \frac{4C_2M}{\nu\varepsilon^2}\nu^{\frac{1}{p}}2^{j\left(-1+\frac{2}{p}+\delta\right)}\big\|\dot\Delta_j \mathcal L^N\big\|_{L^p(0,\mathcal T; L^\infty)}
       &\leq C\left[\frac{\nu 2^{(2+p\delta)j}}{(1-\theta)\frac{\kappa^2}{2}}\right]^{\frac{1}{p}}A~\frac{M\|u_0\|_{\dot{B}_{\infty,\infty}^{-1}}}{\nu\varepsilon^2}~.\nonumber
    \end{align}
    Hence, if $\kappa\neq0$ and $A>1$ satisfy
	\begin{align}\label{kappa condition: multiplicative}
		(1-\theta)\frac{\kappa^2}{2}\geq A^p\left(\frac{M}{\varepsilon^2}\right)^p\left(9C{\bf R}_\infty\right)^p\nu 2^{4N}
	\end{align}
	we would have
	\begin{align}
		\frac{4C_2M}{\nu\varepsilon^2}\nu^{\frac{1}{p}}\big\|\mathcal L^N\big\|_{\widetilde{L}^p\big(0,\mathcal T; \dot{B}_{\infty,\infty}^{-1+2/p}\big)}+\frac{4C_2M}{\nu\varepsilon^2}\nu^{\frac{1}{p}}\big\|\mathcal L^N\big\|_{\widetilde{L}^p\big(0,\mathcal T; \dot{B}_{\infty,\infty}^{-1+2/p+\delta}\big)}\leq 1/9+1/9=2/9~.\nonumber
	\end{align}
    The same thing also holds if one replace $p$ by $p^*$. So in view of (\ref{kappa condition: multiplicative}), we may set $\kappa\neq0$ and $A>1$ such that
    \begin{align}\label{kappa condition: multiplicative'}
       (1-\theta)\frac{\kappa^2}{2}\geq A^{p+p^*}\left(\frac{M}{\varepsilon^2}\right)^{p+p^*}\Big[\left(C{\bf R}_\infty\right)^p+\left(C{\bf R}_\infty\right)^{p^*}\Big]\nu 2^{4N}~,
    \end{align}
    and then $\displaystyle\frac{4C_2M}{\nu\varepsilon^2}\big\|{\mathcal L}^N\big\|_{X^{p,\delta}(\mathcal T)}\leq 4/9$. And so putting together, $\displaystyle\frac{4C_2M}{\nu\varepsilon^2}\|g\|_{X^{p,\delta}(\mathcal T)}\leq 8/9<1$~. Note that to fulfill (\ref{kappa condition: multiplicative'}), we may just set $A=\left(\frac{\kappa^2}{K_\infty^2}\right)^\frac{1}{p+p^*}$ for any pre-given $\kappa$ such that
    \begin{align}
        \kappa^2>K_\infty^2:=\frac{C\nu 2^{4N+1}}{1-\theta}\left(\frac{M}{\varepsilon^2}\right)^{p+p^*}\left({\bf R}_\infty^p+{\bf R}_\infty^{p^*}\right)~.\nonumber
    \end{align}\

   \noindent $\mathbf{\bullet~Proof~of~(\ref{probability lower bound: B^{-1} II}).}$ By the same argument given in the end of Proof of Proposition \ref{Prop. smallness of g II}, one clearly have:
   \begin{align}
       \mathbb P\big(\widetilde\tau\wedge\tau=\infty\big)\geq\left[1-\left(\frac{\kappa}{K_\infty}\right)^{-\frac{2\theta}{p+p^*}}\right]\exp\Bigg(-\frac{6}{\pi^2}\sum_{j>N}j^{-2}\Bigg)~.\nonumber
   \end{align}
   Now if we set $M=e^{\varepsilon_\nu^2T/4}>1$ for $T>0$, then lemma \ref{lem. probability II} says that there exists universal constants $c_0,r_0>0$ such that
   \begin{align}
       \mathbb P\big(\mathcal T_1>T\big)\geq 1-e^{-c_0\varepsilon_\nu^2T}\nonumber
   \end{align}
   whenever $\varepsilon_\nu^2T\geq r_0$. Hence, by independence, we have just proved (\ref{probability lower bound: B^{-1} II}).
\end{proof}\

\subsection{The probability lemma}\label{subsec. The probability lemma}
\begin{proof}[Proof of lemma \ref{lem. probability II}]\
    
    {\noindent$\bf\bullet~Step~I.~Transfer~to~level~set:~$}
    Let $\varepsilon\neq0$, $0<\theta<1$, $M>1$, $T>0$ and the stopping time $\mathcal T_1=\mathcal T_1^{\varepsilon,M}({\rm B})$ be given as in the lemma, we define the set
    \begin{align}
        \Omega_T=\Omega_{T}^{\varepsilon,M}:=\Bigg\{ \sup_{0\leq s\leq t\leq T}\exp{\left[\varepsilon_\nu({\rm B}_t-{\rm B}_s)-\frac{\varepsilon_\nu^2}{4}(t-s)\right]}\leq M\Bigg\}\nonumber
    \end{align}
    so that  $\Omega_T\subset\big\{\mathcal T_1\geq T\big\}$. Next, we fix an arbitrary $0<\alpha<1/2$ and define the $a.s.$-finite random variable
    \begin{align}
        \mathcal H_T^\alpha:=\sup_{0\leq s\leq t\leq T}\frac{{\rm B}_t-{\rm B}_s}{(t-s)^\alpha}~,\nonumber
    \end{align}
    then
    \begin{align}
        \widetilde{\Omega}_T=\widetilde{\Omega}_{T}^{\varepsilon,M}
        :=\Big\{\varepsilon_\nu\mathcal H_T^\alpha(t-s)^\alpha-\frac{\varepsilon_\nu^2}{4}(t-s)\leq\ln M~,\quad \forall~0\leq s\leq t\leq T\Big\}\subset\Omega_T~.\nonumber
    \end{align}
    So if we define the elementary function
    \begin{align}
        \varphi(t)=\varphi^{\varepsilon,M}(t):=\frac{\ln M+\frac{\varepsilon_\nu^2}{4}t}{\varepsilon_\nu t^\alpha}~,\quad 0<t\leq T,\nonumber
    \end{align}
    we see that
    \begin{align}
        \mathcal L_T:=\Big\{\mathcal H_T^\alpha\leq \inf_{0\leq t\leq T}\varphi(t)\Big\}\subset\widetilde{\Omega}_T\subset\Omega_T\subset\big\{\mathcal T_1\geq T\big\}~.\nonumber
    \end{align}
    Hence, the problem transfers to estimating the probability measure of the level set $\mathcal L_T$. By a little work of calculation,
    \begin{align}
        \varphi'(t)
        &=-\frac{\alpha\ln M}{\varepsilon_\nu}t^{-\alpha-1}+\frac{(1-\alpha)\varepsilon_\nu}{4}t^{-\alpha}\nonumber\\
        &=\frac{(1-\alpha)\varepsilon_\nu}{4t^{\alpha+1}}\left(t-\frac{\alpha}{1-\alpha}\frac{\ln M}{\varepsilon_\nu^2/4}\right)~,\nonumber
    \end{align}
    one sees that as long as $t_0:=\frac{\alpha}{1-\alpha}\frac{\ln M}{\varepsilon_\nu^2/4}\leq T$, i.e.
    \begin{align}\label{inf. condition of varphi}
        \ln M\leq \frac{1-\alpha}{\alpha}\cdot\varepsilon_\nu^2T/4~,
    \end{align}
    then
    \begin{align}
        \inf_{0\leq t\leq T}\varphi(t)=\varphi(t_0)=\frac{1}{\alpha}\left(\frac{\alpha}{1-\alpha}\right)^{1-\alpha}\frac{\varepsilon_\nu}{4}\left(\frac{\ln M}{\varepsilon_\nu^2/4}\right)^{1-\alpha}~.\nonumber
    \end{align}
    To fulfill (\ref{inf. condition of varphi}), we may simply set $M$ such that $\ln M=\varepsilon_\nu^2T/4<\frac{1-\alpha}{\alpha}\cdot\varepsilon_\nu^2T/4$ and so
    \begin{align}
        \inf_{0\leq t\leq T}\varphi(t)=\frac{1}{4\alpha^\alpha(1-\alpha)^{1-\alpha}}\cdot\varepsilon_\nu T^{1-\alpha}~.\nonumber
    \end{align}\

    {\noindent$\bf\bullet~Step~II.~Discretisation:~$}
    In this part, we estimate the level set $\big\{\mathcal H_T^\alpha\leq\lambda\big\}$ by showing that
    \begin{align}\label{inclusion for discretisation}
        \bigcap_{n=1}^\infty\bigcap_{k=0}^{2^n-1}\left\{{\rm B}_{(k+1)T/2^n}-{\rm B}_{kT/2^n}\leq\frac{\lambda}{4}\left(\frac{T}{2^n}\right)^\alpha\right\}\subset\big\{\mathcal H_T^\alpha\leq\lambda\big\}~.
    \end{align}
    Assume now
    \begin{align}
        {\rm B}_{(k+1)T/2^n}-{\rm B}_{kT/2^n}\leq D\left(\frac{T}{2^n}\right)^\alpha~,\quad k=0,1,2,...,2^n-1\nonumber
    \end{align}
    holds for all $n\in\mathbb N$. For arbitrarily fixed $0\leq s\leq t\leq T$, we choose $n_0\in\mathbb N$ such that
    \begin{align}
        \frac{T}{2^{n_0}}<t-s\leq 2\cdot\frac{T}{2^{n_0}}~.\nonumber
    \end{align}
    Also, we define the sequences
    \begin{align}
        s_n:=\left\lceil\frac{2^ns}{T}\right\rceil\cdot\frac{T}{2^n}~,\quad t_n:=\left\lfloor\frac{2^nt}{T}\right\rfloor\cdot\frac{T}{2^n}~,\quad n\geq n_0\nonumber
    \end{align}
    i.e. $s_n$ is the smallest dyadic point at level $n$ that is not smaller than $s$; $t_n$ is the largest dyadic point at level $n$ that is not larger than $t$. Note that
    \begin{align}
        0\leq s_n-s_{n+1}\leq\frac{T}{2^{n+1}}~,\quad 0\leq t_{n+1}-t_n\leq\frac{T}{2^{n+1}}~,\nonumber
    \end{align}
    so we obviously have that $s_n\downarrow s$ and $t_n\uparrow t$ as $n\uparrow\infty$. Then we can write
    \begin{align}
        {\rm B}_{s_{n_0}}-{\rm B}_{s}
        &=\sum_{n=n_0}^\infty\left({\rm B}_{s_{n}}-{\rm B}_{s_{n+1}}\right)\nonumber\\
        &\leq D\sum_{n=n_0}^\infty\left(\frac{T}{2^{n+1}}\right)^\alpha\nonumber\\
        &=\frac{D}{1-2^{-\alpha}}\left(\frac{T}{2^{n_0}}\right)^\alpha<D(t-s)^\alpha~.\nonumber
    \end{align}
    Similarly, one can also show that
    \begin{align}
        {\rm B}_{t}-{\rm B}_{t_{n_0}}<D(t-s)^\alpha~.\nonumber
    \end{align}
    On the other hand, by $t_{n_0}-s_{n_0}\leq t-s\leq 2\cdot\frac{T}{2^{n_0}}$, we see that there can be at most one dyadic point at level $n_0$ between $s_{n_0}$ and $t_{n_0}$ and so
    \begin{align}
        {\rm B}_{t_{n_0}}-{\rm B}_{s_{n_0}}\leq 2D\left(\frac{T}{2^{n_0}}\right)^\alpha<2D(t-s)^\alpha~.\nonumber
    \end{align}
    Putting together, we derive
    \begin{align}
        {\rm B}_{t}-{\rm B}_{s}\leq 4D(t-s)^\alpha~.\nonumber
    \end{align}
    Taking $D=\lambda/4$, one proves (\ref{inclusion for discretisation}). Now for each $n$,
    \begin{align}
        \left\{\left(\frac{T}{2^n}\right)^{-1/2}\left({\rm B}_{(k+1)T/2^n}-{\rm B}_{kT/2^n}\right)~:~k=0,1,2,...,2^n-1\right\}\nonumber
    \end{align}
    makes up an independent family of standard normal random variables. Then by the commonly used one-side estimate $\mathbb P\big(N(0,1)>x\big)\leq\frac{1}{2}e^{-x^2/2}~(x\geq0)$, we have that
    \begin{align}
        \mathbb P\Bigg(\bigcup_{k=0}^{2^n-1}\left\{{\rm B}_{(k+1)T/2^n}-{\rm B}_{kT/2^n}>\frac{\lambda}{4}\left(\frac{T}{2^n}\right)^\alpha\right\}\Bigg)\leq\frac{1}{2}2^n\exp\left[-\frac{\lambda^2}{16}\left(\frac{T}{2^n}\right)^{2\alpha-1}\right]~,\quad n\in\mathbb N.\nonumber
    \end{align}
    Hence, (\ref{inclusion for discretisation}) gives that
    \begin{align}\label{probability lower bound of level set H}
        \mathbb P\big(\mathcal H_T^\alpha\leq\lambda\big)\geq 1-\frac{1}{2}\sum_{n=0}^\infty2^n\exp\left[-\frac{\lambda^2}{16}\left(\frac{T}{2^n}\right)^{2\alpha-1}\right]~.
    \end{align}\

    {\noindent$\bf\bullet~Step~III.~Simplify~the~lower~bound:~$}
    To this end, we can set
    \begin{align}
        \lambda=\inf_{0\leq t\leq T}\varphi(t)=\frac{1}{4\alpha^\alpha(1-\alpha)^{1-\alpha}}\cdot\varepsilon_\nu T^{1-\alpha}\nonumber
    \end{align}
    and use the elementary inequality $~2^{n(1-2\alpha)}>1+(2^{1-2\alpha}-1)n>(2^{1-2\alpha}-1)(1+n)\quad (n\geq1)~$ so that
    \begin{align}
        \frac{\lambda^2}{16}\left(\frac{T}{2^n}\right)^{2\alpha-1}
        =\frac{2^{n(1-2\alpha)}}{2^{8}\alpha^{2\alpha}(1-\alpha)^{2-2\alpha}}~\varepsilon_\nu^2T
        >(1+n)c_0\varepsilon_\nu^2T~.\nonumber
    \end{align}
    Here $c_0=c_0(\alpha)=\frac{2^{1-2\alpha}-1}{2^{8}\alpha^{2\alpha}(1-\alpha)^{2-2\alpha}}>0$~. Insert this into (\ref{probability lower bound of level set H}), we see that
    \begin{align}
        \mathbb P\big(\mathcal H_T^\alpha\leq\lambda\big)
        \geq 1-\frac{e^{-c_0\varepsilon_\nu^2T}}{2}\sum_{n=0}^\infty2^ne^{-nc_0\varepsilon_\nu^2T}
        =1-\frac{e^{-c_0\varepsilon_\nu^2T}}{2\big(1-2e^{-c_0\varepsilon_\nu^2T}\big)}~.\nonumber
    \end{align}
    Hence, if take $\varepsilon_\nu^2T$ large enough such that $1-2e^{-c_0\varepsilon_\nu^2T}\geq1/2$~,$~$ i.e. $~\varepsilon_\nu^2T\geq c_0^{-1}\ln4=:r_0$ ,
    \begin{align}
        \mathbb P\big(\mathcal T_1\geq T\big)\geq\mathbb P\big(\mathcal H_T^\alpha\leq\lambda\big)
        \geq 1-e^{-c_0\varepsilon_\nu^2T}~.\nonumber
    \end{align}
\end{proof}

\
\appendix
\renewcommand{\appendixname}{Appendix~\Alph{section}}
\renewcommand{\theequation}{A.\arabic{equation}}

\section{Bernstein type inequality for heat flow and Fourier multipliers behaving like homogeneous functions}
\begin{lemma}[Heat estimates, see \cite{BCD11}]\label{Lem. heat estimates}
	Let $1\leq p\leq p'\leq\infty, 1\leq q\leq q'\leq\infty$, $\lambda>0$ and $0<T\leq\infty$. Let also $\mathcal C'\subset\mathbb R^d$ be some annulus. Then for functions $f\in L^q$ and $g\in L^p(0,T:L^q)$ such that
	\begin{equation}
		{\bf supp}~\widehat f\subset \lambda\mathcal C' \quad and \quad {\bf supp}~\widehat g(t)\subset \lambda\mathcal C'~(0\leq t\leq T),\nonumber
	\end{equation}
	the following estimates hold:
	\begin{align}
        \big\|e^{\nu t\Delta}f\big\|_{L^{q}}
		&\lesssim e^{-ct\lambda^2}\|f\|_{L^q}~\quad (0\leq t\leq T)~,\nonumber\\
		\big\|e^{\nu t\Delta}f\big\|_{L^p(0,T:L^{q'})}
		&\lesssim (\nu\lambda^2)^{-\frac{1}{p}}\lambda^{\frac{d}{q}-\frac{d}{q'}}\|f\|_{L^q}~,\nonumber\\
		\left\|\int_0^t e^{\nu(t-s)\Delta}g(s){\rm d}s\right\|_{L^{p'}(0,T:L^{q'})}
		&\lesssim (\nu\lambda^2)^{-1+\frac{1}{p}-\frac{1}{p'}}\lambda^{\frac{d}{q}-\frac{d}{q'}}\|g\|_{L^p(0,T:L^{q})}~.\nonumber
	\end{align}
	Here, constant $c>0$ depends only on the annulus $\mathcal C'$ while the implicit constants can depend on the annulus $\mathcal C'$ and dimension $d$.
\end{lemma}

\begin{lemma}[See \cite{BCD11}] \label{Lem. homogeneous multipliers}
    Let $\mathcal C$ be an annulus, $m\in\mathbb R$ and $k=2\big\lfloor 1+d/2\big\rfloor$. Let $\sigma$ be a $k$-times differentiable function on $\mathbb R^d\setminus\{0\}$ such that: for any multi-index $|\alpha|\leq k$ there exists a constant $c_\alpha>0$ such that
    \begin{align}
        |\partial^\alpha\sigma(\xi)|\leq c_\alpha|\xi|^{m-|\alpha|}~,\quad \xi\in\mathbb R^d\setminus\{0\}~.\nonumber
    \end{align}
    Then there exists a constant $C=C(c_\alpha)>0$ such that: for any $q\in[1,\infty]$, $\lambda>0$ and any function $f\in L^q$ $s.t.~{\bf supp}~\widehat f\subset \lambda\mathcal C$, we have
    \begin{align}
        \|\sigma(D)f\|_{L^q}\leq C\lambda^m\|f\|_{L^q}~.\nonumber
    \end{align}
    Here $\sigma(D)f:=\mathcal F^{-1}\big(\sigma \hat{u}\big)$.
\end{lemma}

\bibliographystyle{alpha}
\bibliography{Paper2ref}

@book{BCD11,
	author = {Hajer Bahouri and Jean-Yves Chemin and Rapha\"el Danchin},
	publisher = {Springer},
	title = {{Fourier Analysis and Nonlinear Partial Differential Equations}},
	volume = {343},
	year = {2011}}

@article{CG09,
	author = {J.-Y. Chemin and I. Gallagher},
	journal = {Ann. I. H. Poincar{\'e}---- AN},
	number = {2009},
	pages = {599--624},
	title = {{Well-posedness and stability results for the Navier-Stokes equtions in $\mathbb R^3$}},
	volume = {26},
	year = {2009}}

@book{CMP93,
	author = {M. Cannone and Y. Meyer and F. Planchon},
	publisher = {{S\'eminaire ''\'Equations aux D\'eriv\'ees Partielles'' de l'\'Ecole polytechnique, Expos\'e VIII}},
	title = {Solutions autosimilaires des\'equations de Navier-Stokes},
	year = {1993-1994}}

@article{BP08,
	author = {J. Bourgain and N. Pavlovi\'c},
	journal = {Journal of Functional Analysis},
	number = {2008},
	pages = {2233-2247},
	title = {{Ill-posedness of the Navier-Stokes equations in a critical space in 3D}},
	volume = {255},
	year = {2008}}

@article{KT01,
	author = {H. Koch and D. Tataru},
	journal = {Adv. Math.},
	pages = {22--35},
	title = {{Well-posedness for the Navier-Stokes equations}},
	volume = {157},
	year = {2001}}

@book{Can03,
	author = {M. Cannone},
	publisher = {Friedlander, S. \& Serre, D. (eds.)},
	series = {Handbook of mathematical fluid dynamics},
	title = {{Harmonic analysis tools for solving the incompressible Navier-Stokes equations}},
	year = {2003}}

@article{Che92,
	author = {J.-Y. Chemin},
	journal = {SIAM J. Math. Anal.},
	pages = {20--28},
	title = {{Remarques sur l'existence globale pour le syst{\`e}me de Navier-Stokes incompressible}},
	volume = {23},
	year = {1992}}

@article{FK64,
	author = {H. Fujita and T. Kato},
	journal = {I. Arch. Ration. Mech. Anal.},
	pages = {269--315},
	title = {{On the Navier--Stokes initial value problem}},
	volume = {16},
	year = {1964}}

@article{FR02,
	author = {Franco Flandoli and Marco Romito},
	journal = {Trans. Amer. Math. Soc.},
	month = {02},
	number = {6},
	pages = {2207-2241},
	title = {{Partial regularity for the stochastic Navier-Stokes equations}},
	volume = {354},
	year = {2002}}

@article{GHV14,
	author = {Nathan E. {Glatt}-{Holtz} and Vlad Vicol},
	journal = {The Annals of Probability},
	month = {01},
	number = {1},
	pages = {80-145},
	title = {{Local and global existence of smooth solutions for the stochastic Euler equations with multiplicative noise}},
	volume = {42},
	year = {2014}}

@article{RZZ14,
	author = {Michael R\"ockner and Rongchan Zhu and Xiangchan Zhu},
	journal = {Stochastic Process. Appl.},
	number = {5},
	pages = {1974-2002},
	title = {{Local existence and non-explosion of solutions for stochastic frational partial differential equations driven by multiplicative noise}},
	volume = {124},
	year = {2014}}

@article{FL21,
	author = {Franco Flandoli and Dejun Luo},
	journal = {Probability Theory and Related Fields},
	month = {03},
	number = {1-2},
	pages = {309-363},
	title = {{High mode transport noise improves vorticity blow-up control in 3D Navier-Stokes equations}},
	volume = {180},
	year = {2021}}

@article{DPD03,
	author = {Giuseppe {Da Prato} and A. Debussche},
	journal = {J. Math. Pures Appl.},
	pages = {877-947},
	title = {{Ergodicity for 3D stochastic Navier-Stokes equations}},
	volume = {82},
	year = {2003}}

@article{ABC21,
	author = {Dallas Albritton and Elia Bru{\'e} and Maria Colombo},
	journal = {Ann. of Math.},
	month = {May},
	number = {1},
	pages = {415-455},
	title = {{Non-uniqueness of Leray solutions of the forced Navier-Stokes equations}},
	volume = {196},
	year = {2022}}

@article{BCV18,
	author = {Tristan Buckmaster and Maria Colombo and Vlad Vicol},
	journal = {J. Eur. Math. Soc.},
	month = {Oct.},
	number = {9},
	pages = {3333-3378},
	title = {{Wild solutions of the Navier-Stokes equations whose singular sets in time have Hausdorff dimension strictly less than 1}},
	volume = {24},
	year = {2022}}

@article{FR08,
	author = {Franco Flandoli and Marco Romito},
	doi = {10.1007/s00440-007-0069-y},
	journal = {Probability Theory and Related Fields},
	number = {3-4},
	pages = {407--458},
	publisher = {Springer Science and Business Media {LLC}},
	title = {{Markov selections for the 3D stochastic Navier{\textendash}Stokes equations}},
	url = {https://doi.org/10.1007%2Fs00440-007-0069-y},
	volume = {140},
	year = 2008}

@article{Kat84,
	author = {Tosio Kato},
	journal = {Math. Z.},
	number = {4},
	pages = {471-480},
	title = {{Strong ${L}^p$-solutions of the Navier-Stokes equation in $\bf{R}^m$, with applications to weak solutions}},
	volume = {187},
	year = {1984}}

@article{BV19,
	author = {Tristan Buckmaster and Vlad Vicol},
	eprint = {1901.09023},
	journal = {EMS Surv. Math. Sci.},
	number = {1/2},
	pages = {173-263},
	title = {{Convex integration and phenomenologies in turbulence}},
	url = {https://arxiv.org/pdf/1901.09023.pdf},
	volume = {6},
	year = {2019}}

@article{CL22,
	author = {Alexey Cheskidov and Xiaoyutao Luo},
	journal = {Invent. Math.},
	month = {05},
	pages = {987-1054},
	title = {{Sharp nonuniqueness for the Navier-Stokes equations}},
	volume = {229},
	year = {2022}}

@article{BVannals19,
	author = {Tristan Buckmaster and Vlad Vicol},
	journal = {Ann. of Math.},
	number = {1},
	pages = {101-144},
	title = {{Nonuniqueness of weak solutions to the Navier-Stokes equation}},
	volume = {189},
	year = {2019}}

@article{ABC23,
	author = {D. Albritton and E. Brué and M. Colombo},
	journal = {Ann. PDE},
	number = {17},
	title = {{Gluing Non-unique Navier–Stokes Solutions}},
	volume = {9},
	year = {2023}}

@article{CP25,
  author = {M. P. Coiculescu and S. Palasek},
  title = {Non-uniqueness of smooth solutions of the Navier–Stokes equations from critical data},
  journal = {Inventiones mathematicae},
  year = {2025},
  volume = {244},
  issn = {1432-1297},
  url = {http://dx.doi.org/10.1007/s00222-025-01396-z},
  doi = {10.1007/s00222-025-01396-z},
  number = {1},
  publisher = {Springer Science and Business Media LLC},
  month = Dec,
  pages = {165–219},
}

@article{Ver81,
  author = {A. J. Veretennikov},
  title = {On Strong Solutions and Explicit Formulas for Solutions of Stochastic Integral Equations},
  journal = {Mathematics of the USSR-Sbornik},
  year = {1981},
  volume = {39},
  issn = {0025-5734},
  url = {http://dx.doi.org/10.1070/SM1981v039n03ABEH001522},
  doi = {10.1070/sm1981v039n03abeh001522},
  number = {3},
  publisher = {Steklov Mathematical Institute},
  month = Apr,
  pages = {387–403},
}

@misc{Agresti2026,
  author = {Antonio Agresti},
  title = {{On the absence of blow-up in the 3D Navier-Stokes equations with transport noise}},
  year = {2026},
  doi = {10.48550/ARXIV.2607.15140},
  url = {https://arxiv.org/abs/2607.15140},
  publisher = {arXiv},
  copyright = {arXiv.org perpetual, non-exclusive license},
}

@misc{HLL2025,
  author = {W. Hong and S. Li and W. Liu},
  title = {{Stochastic Forced 3D Navier-Stokes Equations in $\mathbb{H}^{1/2}$-Space}},
  year = {2025},
  doi = {10.48550/ARXIV.2511.15223},
  url = {https://arxiv.org/abs/2511.15223},
  publisher = {arXiv},
  copyright = {arXiv.org perpetual, non-exclusive license},
}

\end{document}